\documentclass{article}
\usepackage[margin=30mm]{geometry}
\usepackage{color}
\usepackage{graphicx} 
\usepackage{amsthm, amsmath,color,amssymb}
\usepackage{url}
\usepackage{algorithm}
\usepackage{algpseudocode}
\usepackage{capt-of}

\usepackage{tikz}
\usetikzlibrary{intersections,calc,arrows.meta,patterns}

\title{Comprehensive Restriction Algorithm for Hypergeometric Systems}
\author{Hiromasa Nakayama and Nobuki Takayama}
\date{September 24, 2026}

\def\aa{a}
\def\bb{b}
\def\cc{c}

\def\CC{\mathbb{C}}
\def\KK{{\bf K}}
\def\NN{\mathbb{N}}
\def\QQ{\mathbb{Q}}
\def\RR{\mathbb{R}}
\def\ZZ{\mathbb{Z}}
\def\memo#1{{\color{blue} #1}}
\def\p{\partial}
\def\pd#1{\partial_{#1}}
\def\qed{\hfill $\Box$ \medbreak}

\def\CGS{{\rm CGS}}

\def\aalpha{a}
\def\bbeta{b}
\def\ggamma{c}
\def\Hg{H_g}

\newtheorem{theorem}{Theorem}
\newtheorem{proposition}{Proposition}
\newtheorem{example}{Example}
\newtheorem{remark}{Remark}

\begin{document}

\maketitle

\section{Introduction}

We denote by $D$ the Weyl algebra
$$\CC \langle x_1, \ldots, x_n, \pd{1}, \ldots, \pd{n} \rangle, $$
that is the ring of linear partial differential operators
with polynomial coefficients.
When we need to specify the number of variables, we denote $D$ by $D_n$.
Let $M$ be a holonomic $D$-module 
on the $n$-dimensional space $\CC^n=\{ x=(x_1, \ldots, x_n) \}$.
The $0$-th restriction of $M$ to $V(x_{m+1}, \ldots, x_n)$ 
is defined as
$$
\frac{M}{x_{m+1}M + \cdots + x_n M}
$$
(see, e.g., \cite{bjork-2012}, \cite[Chap 5]{SST}).
An algorithm computing the restriction 
was given by T.Oaku \cite{oaku-1997}.
In this paper, we consider a problem of computing the restriction
for a given holonomic $D$-module with parameters.
We will give a partial answer to the problem for general holonomic $D$-modules
and an answer to hypergeometric holonomic $D$-modules.

The basic method for performing various calculations on ideals 
or submodules of free modules
involving parameters is 
the comprehensive Gr\"obner basis introduced 
by V.Weispfenning \cite{weispfenning-1992}.
K.Nabeshima, K.Ohara, S.Tajima \cite{not-2016} introduced comprehesive Gr\"obner systems
(CGS) for rings of linear partial differential operators.
They applied their method of computing CGS to the problem of computing
$b$-functions for polynomials with parameters.
The parameter space is stratified so that a $b$-function is associated
to each stratum.

For a given holonomic $D$-module with parameters,
we want to stratify the parameter space so that
a restriction module that does not depend on parameters
is associated to each stratum.
We start with generalizing the method by K.Nabeshita et al. 
to compute a generic $b$-function
that is also called an indicial polynomial or 
a $b$-function for restriction (Section \ref{sec:comprehensive}).
The maximal integral root of it plays the central role
to apply the Oaku's $b$-function criterion 
of the restriction algorithm \cite{oaku-1997}.
The parameter space can be stratified so that a generic $b$-function
is associated to each stratum.
However, the difficulty is that roots of it still depends on parameters.
We use isomorphic correspondences of $D$-modules with parameters
to address this difficulty.
Since we only need to consider integral roots to obtain the restriction, 
isomorphic correspondences are fully available.

We focus on algorithms to construct isomorphisms
among hypergeometric $D$-modules in sections \ref{sec:contiguity1},
\ref{sec:rat-sol}, \ref{sec:restriction_of_A_hg}, \ref{sec:isom}.
M.Saito gave a method to classify GKZ hypergeometric systems
into isomorphic classes \cite{saito-isom}.
By making this method a complete algorithm and providing a restricting algorithm for isomorphisms,
we give an algorithm to classify a class of hypergeometric systems of Horn type
\cite{berkesch} into isomorphic classes.
The key ingradient of our method is constructing strata so that
a contiguity relation (a left $D$-isomorphism) of a hypergeometric system with parameters
is associated to each stratum.
Note that a general algorithm to check if two holonomic $D$-modules
are isomorphic or not is given by H.Tsai and U.Walther \cite{tw-hom}.
Considering a comprehensive version of this algorithm is a future problem.

Utilizing algorithms to classifying isomorphic classes
of hypergeometric systems, we finally give a comprehensive restriction
algorithm in section \ref{sec:restriction}.
The remaining sections \ref{sec:gauss-hg}, \ref{sec:f1} are discussions 
on restrictions to the origin
of the Gauss hypergeometric system and the Appell $F_1$ system.

\section{Comprehensive Gr\"obner System and Generic $b$-function} \label{sec:comprehensive}

K.Nabeshima, K.Ohara, S.Tajima introduced
an algorithm for computing comprehensive Gr\"obner systems (CGS) in rings of linear partial differential operators \cite{not-2016}.
They also gave applications of CGS for computing $b$-functions for singularities.
We apply their algorithm to obtain $b$-functions 
for weight vectors
to compute restrictions of $D$-modules. See, e.g., \cite[Chap 5]{SST} on $b$-functions for weight vectors. Being inspired by the computer algebra system 
Risa/Asir \cite{risa-asir} command name
{\tt generic\_bfct}, we call them {\it generic $b$-functions}.
We also call a generic $b$-function {\it a $b$-function for restriction}
in this paper to distinguish with a $b$-function of a contiguity relation
and a $b$-function for a polynomial.

Let 
$$D_n[\beta]=\CC[\beta_1, \ldots, \beta_m]
  \langle x_1, \ldots, x_n, \pd{1}, \ldots, \pd{n} 
  \rangle
$$
be the Weyl algebra with parameters
$\beta=(\beta_1, \ldots, \beta_m)$
regarded as indeterminates.
We denote $D_n$ by $D$ when the number of variables is clear.
For a left ideal $I$ generated by a set
of generators $P$ in $D[\beta]$,
we compute a Gr\"obner basis $G$ with a block order
$\succ_b$
satisfying $x_i, \pd{i} \succ_b \beta_j$
for any $i$ and $j$
where $\succ$ is a tie-breaker of the block order.
Put $E=G \cap \CC[\beta]$ and $G'=G \setminus E$.
The set $V(E)$ is a given set of equality constraints on parameters.
We denote by $\CGS(E, N, P, \succ_b)$
or by $\CGS(E, N, I, \succ_b)$
the output CGS of $I$ on $V(E) \setminus V(N)$
where $N$ is a given set of equality constraints on parameters.
The CGS is a finite set of data of the form $(V(E_i) \setminus V(N_i), {\cal G}_i)$
where $E_i, N_i \subset \CC[\beta]$, ${\cal G}_i \subset D[\beta]$ and
they are finite set.
The CGS has a property that
for any $a \in V(E_i) \setminus V(N_i)$,
${\cal G}_i |_{\beta=a}$ is a Gr\"obner basis of $P|_{\beta=a}$ in $D$ 
with respect to the order $\prec$.
$V(E_i) \setminus V(N_i)$ is called a stratum 
and the strata of this form in the CGS cover $V(E) \setminus V(N)$.
Note that when $E=\emptyset$, we regard $V(E)=\CC^m$.
The procedure $\CGS$ is recursively called.
At the top level, we usually start with $\CGS(E=\emptyset,N=\{1\},P,\succ_b)$.
Note that $V(\{1\})=\emptyset$.


Let $w$ be a vector in $\ZZ_{\geq 0}^n$.
The $(-w,w)$-degree of $a x^p \pd{}^q$ 
is $-w\cdot p + w \cdot q$
where $a \in \CC[\beta_1, \ldots, \beta_m]$,
$a \not= 0$, and 
$x^p=\prod_{i=1}^n x_i^{p_i}$, 
$\pd{}^q=\prod_{i=1}^n \pd{i}^{q_i}$. 
The $(-w,w)$-initial term for 
$\ell=\sum_{(p,q)\in E} a_{pq} x^p \pd{}^q$
is the sum of the maximal $(-w,w)$-degree terms of $\ell$
and is denoted by 
${\rm in}_{(-w,w)}(\ell)$.
For a given left ideal $I$ in $D[\beta]$,
the ideal generated by ${\rm in}_{(-w,w)}(\ell)$, $\ell \in I$
is called {\it the initial form ideal} (with respect to the weight
vector $(-w,w)$).
The {\it $(-w,w)$ order} of $\ell$ is defined by 
${\rm max}_{(p,q) \in E}\, -wp+wq$
and is denoted by
${\rm ord}_{(-w,w)}(\ell)$.

Let $\succ$ be a term order in $D[\beta]$.
The order $\succ_{(-w,w)}$ is defined as
\begin{eqnarray*}
&&  c_{pq}(\beta) x^p \pd{}^q \succ_{(-w,w)}  c'_{p'q'}(\beta) x^{p'} \pd{}^{q'} \\
&\Leftrightarrow&
 -w\cdot p + w \cdot q >  -w\cdot p' + w \cdot q' \\
 & & \quad \mbox{ or } \ 
 (-w\cdot p + w \cdot q =  -w\cdot p' + w \cdot q' \ \mbox{ and } \ 
  x^p \pd{}^q \succ x^{p'} \pd{}^{q'}). 
\end{eqnarray*}
Since the order $\succ_{(-w,w)}$ is not a well-order,
we need to utilize the homogenized Weyl algebra to compute
Gr\"obner bases with this order.
Their CGS algorithm can also be applied 
to the homogenized Weyl algebra with parameters
(see, e.g., \cite[Th. 1.2.6]{SST} on the homogenized Weyl algebra),
and obtain a CGS for the initial form ideal
of a given left ideal in $D[\beta]$.
This method utilizing the homogenized Weyl algebra 
is not explicitly described in the paper \cite{not-2016}, 
so we explain it below.
Note that the case of the holonomic $D$-module $M$ is of the form
$D^m/I$ where $I$ is a submodule can be discussed analogously.

\begin{algorithm}[H]\rm 
\footnotesize
\caption{Computing parametric initial form ideal} \label{alg:CGS}
~~
\begin{itemize}
\item {\bf Require}:  a set of generators of a left ideal $I$
in $D[\beta]$,~ a weight vector $w \in (\ZZ_{\geq 0})^n$
\item {\bf Ensure}:
  A stratification of the parameter space $\{(E_i,N_i)\}$ and
generators of the initial form ideal ${\rm in}_{(-w,w)}(I)$ on the stratum $V(E_i) \setminus V(N_i)$.
\end{itemize}
\begin{enumerate}

\item 
Let $\succ^h_{(-w,w)}$ be an order in the homogenized Weyl algebra
defined as 
\begin{eqnarray*}
&&  c_{pq}(\beta) x^p \pd{}^q h^r \succ^h_{(-w,w)}  c'_{p'q'}(\beta) x^{p'} \pd{}^{q'} h^{r'} \\
&\Leftrightarrow&
 -w\cdot p + w \cdot q >  -w\cdot p' + w \cdot q' \\
 & & \quad \mbox{ or } \ 
 (-w\cdot p + w \cdot q =  -w\cdot p' + w \cdot q' \ \mbox{ and } \ 
  x^p \pd{}^q h^r \succ x^{p'} \pd{}^{q'} h^{r'}). 
\end{eqnarray*}
where the tie-breaker is an elimination order of $h$.
Extending $\succ^h_{(-w,w)}$ to a block order $\succ^h_{b,(-w,w)}$
such that $x_i, \pd{i}, h \succ^h_{b,(-w,w)} \beta_j$,
we compute a CGS 
$$
\mathcal{G} = \{(E_i, N_i, {\cal G}_i) \mid i = 1, 2, \ldots, m\}
$$
for $I^h$.

\item 
Return $\{ ({\cal G}_i|_{h=1}, E_i, N_i)\,|\, i=1, \ldots, m\}$.
\end{enumerate}
\normalsize
\end{algorithm}
Polynomials in only $h$ and $\beta$ obtained 
by eliminating $x_1, \ldots, \pd{n}$ and setting $h=1$
define varieties $N_i$.
Note that the homogenization variable $h$ appears only in this section. 
In subsequent sections, the symbol $h$ will be used for other purposes, such as a element of a Hilbert basis.

\begin{example} \rm  \label{ex:f1-f2-w-w-gb}
Our first example is
the system of Appell differential operators for $F_1(a,b,b',c;x,y)$,
that is
\begin{eqnarray}
&&x(1-x)\pd{x}^2+y(1-x)\pd{x}\pd{y}+(\cc-(\aa+\bb+1)x)\pd{x}-\bb y \pd{y} - \aa \bb,  \label{eq:f1a}\\
&&y(1-y)\pd{y}^2+x(1-y)\pd{x}\pd{y}+(\cc-(\aa+\bb'+1)y)\pd{y}-\bb' x \pd{x} - \aa \bb',  \label{eq:f1b} \\
&&(x-y) \pd{x}\pd{y} - \bb' \pd{x} + \bb \pd{y}. \label{eq:f1c}
\end{eqnarray}
They annihilate the function $F_1$.
The left ideal $I$ generated by them are holonomic ideal for any value
of the parameter vector and $D_2/I$ is a holonomic $D_2$-module for any specialization of the parameter vector.
For the weight vector $(-w,w), w = (1,1)$,
we apply Algorithm \ref{alg:CGS} to obtain CGS with respect to 
the order $\succ_{(-w,w)}$.
The parametric initial form ideal
is generated for any parameter values $(a,b,b',c)$
by 
\begin{align*}
&(x-y) \pd{x} \pd{y} + b \pd{y} - b' \pd{x}, \\
&-y \pd{x}\pd{y} - y \pd{y}^2 - b' \pd{x} +(b-c) \pd{y}, \\
&-x \pd{x}^2 + y \pd{y}^2 +(b'-c) \pd{x} + (-b+c) \pd{y}, \\
& (-xy+y^2)\pd{y}^2 - b' x \pd{x} + (b-c)x\pd{y} + cy \pd{y}.   
\end{align*}

Our second example is
the system of Appell differential operators for
$F_2(a,b,b',c,c';x,y)$, that is
\begin{eqnarray}
&&x(1-x)\pd{x}^2-xy\pd{x}\pd{y}+(\cc-(\aa+\bb+1)x)\pd{x}-\bb y \pd{y}-\aa \bb, \label{eq:f2a} \\
&&y(1-y)\pd{y}^2-xy\pd{x}\pd{y}+(\cc'-(\aa+\bb'+1)y)\pd{y}-\bb' x \pd{x}-\aa \bb'. \label{eq:f2b}
\end{eqnarray}
For the weight vector $(-w,w), w = (1,1)$,
parametric initial form ideal
is generated for any parameter values $(a,b,b',c)$
by 
\begin{align*}
& -x^2 y \pd{x}^3 \pd{y} + xy(x-y) \pd{x}^2 \pd{y}^2 + x y^2 \pd{x}\pd{y}^3 
-b' x^2 \pd{x}^3 + (c'x - (a+b'+c+3)y) x \pd{x}^2 \pd{y} + \\
& ((a+b+c'+3)x - cy)x\pd{x} \pd{y}^2
+ b y^2 \pd{y}^3 - (a + c + 2) b' x \pd{x}^2 + ((a+b+2)c'x-(a+b'+2)cy)\pd{x}\pd{y} + \\
& (a+c'+2)by \pd{y}^2 - (a+1)b'c \pd{x} + (a+1)bc' \pd{y}, \\
&y \pd{y}^2 + c' \pd{y}, \\
&x \pd{x}^2 + c \pd{x}. 
\end{align*}
In the above examples, there is only one stratum.

The third example is the left ideal generated by
$a x \pd{x} + b y \pd{y}$ and $x \pd{x}+y\pd{y}$
where $a, b$ are parameters.
When $a-b \not= 0$, the $(-1,-1,1,1)$ initial form ideal is generated
by $x \pd{x}$ and $y \pd{y}$.
When $a-b=0$, it is generated by $x \pd{x}+y\pd{y}$.
There are two strata.
\end{example}

The algorithm \cite[Th.5.1.6.]{SST} for computing the generic $b$-function for any weight vector of a holonomic $D$-ideal can be generalized to ideals with parameters in coefficients as follows.

\begin{algorithm}[H]\rm 
\footnotesize
\caption{Parametric generic $b$-function} \label{alg:CGS-bf}
~ 
\begin{itemize}
\item 
{\bf Require}:  A set of generators $P$ 
of a holonomic left ideal $I$ in $D[\beta]$,~ a weight vector $w \in \ZZ_{\geq 0}^n$
\item
{\bf Ensure}: Stratification and the generic $b$-function $b(s)$  on each stratum where
$\langle b(s) \rangle = {\rm in}_{(-w,w)}(I) \cap \CC[s]$ ($s = \sum_{i=1}^n w_i \theta_i$, $\theta_i=x_i \pd{i}$). 
\end{itemize}
\begin{enumerate}
\item 
Compute a parametric initial form ideal ${\rm in}_{(-w,w)}(I)$.
We obtain
the initial form ideal generated by $\mathcal{G}_i$
on each stratum $(E_i, N_i)$ $(i = 1, \ldots, r)$.
\item 
$B \leftarrow \emptyset$ 
\item 
For each $i = 1, \ldots, r$, do
	\begin{enumerate}
	\item[3.1] 
	For each element $\ell$ of $\mathcal{G}_i$,
	make a replacement $x_k \rightarrow u_k x_k$, $\pd{k} \rightarrow v_k \pd{k}$
        where $k$ runs over a set of indices such that $w_k \neq 0$.
        Let $J_i$ be the left ideal geneted by these $\ell$'s and $1 - u_k v_k$.
	\item[3.2] 
	Compute a CGS for the left ideal $J_i$ on the stratum $(E_i, N_i)$.
         We use an elimination order $\succ$ of $u_k, v_k$.
         (call $\CGS(E_i,N_i,J_i, \succ_b)$.)
	Let $\mathcal{G}_{ij}$ and stratum $(E_{ij}, N_{ij})$ $(j = 1, \ldots, n_i)$ be the output.
    Collect all elements that do not contain $u_k, v_k$
    from $\mathcal{G}_{ij}$ and put them in $\mathcal{G}_{ij}'$. 
	\item[4] 
	For each $j = 1, \ldots, n_i$, do
		\begin{enumerate}
		\item[4.1] 
		Any element $P$ of $\mathcal{G}_{ij}'$ is of the form $P = x^a p(\theta) \p^b$.
		 Replace it as $[\theta]^a p(\theta - b) [\theta]_b$ and put $J_{ij}$ the ideal generated by them where $\theta_i = x_i \pd{i}$.
 Here, $[\theta]^a=\prod_{j=1}^n \prod_{l=1}^{a_{j}} (\theta_j+l)$
and  $[\theta]_b=\prod_{j=1}^n \prod_{l=0}^{b_{j}-1} (\theta_j-l)$,
\cite[p.45, p.195]{SST}.
		\item[4.2] 
		Add $s - \sum_{i=1}^n w_i x_i \p_i$ to the ideal $J_{ij}$.
                Regard it as an ideal in 
		$\CC[\beta]\langle \theta_1, \ldots, \theta_n, s \rangle$,
                compute a CGS, 
                and obtain the generator of 
                $J_{ij} \cap \mathbb{C}[s]$ on the
		stratum $(E_{ijk}, N_{ijk})$.
		In other words, compute $\CGS(E_{ij}, N_{ij}, J_{ij}, \succ'_b)$ where  $\succ'$
is an order satisfying $x, \pd{x} \succ' s$
               and take the minimal degree polynomial $b(s)$ of $s$ with coefficients
               in $\CC[\beta]$ for each stratum. 
               Add the polynomial $b(s)$ and the stratum to $B$.
		\end{enumerate}
	\end{enumerate}
Return $B$.
\end{enumerate}
\normalsize
\end{algorithm}

\begin{example} \rm  \label{ex:f1-f2-comprehensive-b-fn}
The generic $b$-function for $(-w,w)$, $w=(1,1)$
of the Appell system of $F_1(a,b,b',c)$ 
is
$$
b(s) = s(s+c-1)
$$
on $\CC^4=\{(a,b,b',c)\}$.

The generic $b$-functions for $(-w,w)$, $w=(1,1)$
of the Appell system of $F_2(a,b,b',c,c')$ are
\begin{align*}
& \mbox{stratum} & & \mbox{ generic $b$-function} \\
& V(0) \setminus V((c-c')(c+c'-2)) & &  s(s+c-1)(s+c'-1)(s+c+c'-2) \\
& V(c-c') & & s(s+c'-1)(s+2c'-2) \\
& V(c+c'-2) \setminus V(c-c') & & s(s-c'+1)(s+c'-1) 
\end{align*}
\end{example}

\section{Review on Algorithms for Contiguity Relations} \label{sec:contiguity1}

In this section, we review known algorithms to find contiguity relations.
In the sections \ref{sec:rat-sol} and \ref{sec:restriction_of_A_hg},
we propose new algorithms to find contiguity relations.

Let $\KK$ be a rational function field
$\CC(\beta_1, \ldots, \beta_d)$.
Let 
$D_n=\KK \langle x_1, \ldots, x_n, \pd{1}, \ldots, \pd{n} \rangle$ 
be the Weyl algebra of $n$ variables over the field $\KK$.
We denote $D_n$  by $D$ when the number of variables is clear. 
We consider a family of holonomic $D$-modules $M(\beta)=D/I(\beta)$
where $I(\beta)$ is a left ideal of $D$.
The parameters $\beta_i$'s are specialized to complex numbers
in some context.

In what follows in this section, we assume that the parameters are
specialized to be numbers.
Let $H_i(\beta)$ be an element of $D$
satisfying the condition
\begin{equation}  \label{eq:H-relation}
 \ell \, H_i(\beta) \in I(\beta) \mbox{ for all $\ell \in I(\beta+e_i)$}.
\end{equation}
Here $\beta+e_i$ means $(\beta_1, \ldots, \beta_{i-1},\beta_i+1,\beta_{i+1}, \ldots, \beta_d)$.
Then, we have the left $D$-morphism
\begin{equation} \label{eq:d-morphism}
 \varphi_i : M(\beta+e_i)  \ni [p] \longrightarrow [p H_i(\beta)] \in M(\beta).
\end{equation}
The morphism $\varphi_i$ induces the morphism of vector spaces of the opposite direction
\begin{equation} \label{eq:opp-d-morphism}
 {\rm Hom}_D(M(\beta),{\hat {\cal O}}_a) \ni f \longrightarrow H_i(\beta)  f \in {\rm Hom}_D(M(\beta+e_i),{\hat {\cal O}}_a)
\end{equation}
where ${\hat {\cal O}}_a$ is a germ of formal power series at a point $x=a$.
When the morphism (\ref{eq:d-morphism}) is an isomorphism,
the opposite linear map (\ref{eq:opp-d-morphism}) is also an isomorphism.
The operator $H_i(\beta)$ is called the {\it up-step operator} (for the direction $i$)
or the up-step contiguity operator.
Analogously, if we have an element $B_i \in D$
satisfying
\begin{equation} \label{eq:B-relation}
 \ell \, B_i(\beta+e_i) \in I(\beta+e_i) \mbox{ for all $\ell \in I(\beta)$},
\end{equation}
then we have a left $D$-morphism
\begin{equation} \label{eq:d-morphism-down}
 \psi : M(\beta)  \ni [p] \longrightarrow [p B_i(\beta+e_i)] \in M(\beta+e_i),
\end{equation}
the operator $B_i(\beta)$ is called the {\it down-step operator} 
or the down-step contiguity operator\footnote{The symbols H and B are a convention in the theory of special functions, derived from the French terms ``en haut'' and ``en bas''.}.

Regard $\beta$ as indeterminates. 
We consider the composite
\begin{equation}  \label{eq:b-function-for-contiguity}
  M(\beta) \stackrel{B_i(\beta+e_i) H_i(\beta)}{\longrightarrow} M(\beta),
\quad
  M(\beta+e_i) \stackrel{H_i(\beta)B_i(\beta+e_i)}{\longrightarrow} M(\beta+e_i)
\end{equation}
When they are multiplications of a polynomial in $\beta$, 
it is called {\it a $b$-function of the contiguity} among $\beta$ and $\beta+e_i$.
When the value of the $b$-function is not zero at a value of $\beta$,
contiguity operators give an isomorphism among
$M(\beta)$ and $M(\beta+e_i)$.
We call the set of up-step operator, down-step operator, and the $b$-function 
{\it contiguity relation}\/.

We note that the same name of $b$-function is also used in the previous section
in a different context.
If there is a risk of confusions,
we call the $b$-function in the previous section
{\it the $b$-function for restriction} and the $b$-function in this section
{\it the $b$-function of the contiguity}.
Note that the letter $b$ is sometimes used to represent a parameter, in keeping with the traditional convention for hypergeometric functions; the following discussion may also use $b$ as a parameter.

\begin{example}\rm
We consider the Gauss hypergeometric operator
with $a=c, b, c$
and denote it by
\begin{equation}
L(b,c)=x(1-x)\pd{x}+(c-(c+b-1)x)\pd{x}-c b.
\end{equation}
Put $\beta_1=b, \beta_2=c$, 
$I(\beta)=DL(b,c)$ and
consider $M(\beta)=D/I(\beta)$.
We fix
$b$ as a generic complex number.
We will derive a contiguity relation with respect to $c$.
Put $\theta_x = x \pd{x}$. 
The operator $\theta_x$ is called the Euler operator.
Since
$$ x L(c) = \theta_x(\theta_x+c-1)-x (\theta_x+c)(\theta_x+b)
  = (\theta_x+c-1) (\theta_x - x (\theta_x+b)),
$$
the classical solution space of it is spanned by
\begin{eqnarray*}
f_1(c) &=& (1-x)^{-b} \\
f_2(c) &=& x^{1-c} {}_2F_1(1,1+b-c,2-c;x)
\end{eqnarray*}
as a vector space over $\CC$.
An up-step operator and a down-step operator with respect to $c$ are
\begin{eqnarray}
H(c) &=& (x-1)\pd{x}+c \\
B(c) &=& (1-c)\left( x (x-1) \pd{x} + b x - c + 1\right). 
\end{eqnarray}
The $b$-function for the contiguity is 
\begin{equation}
c^2 (c-b). 
\end{equation}
These operators act to solutions as follows.
\begin{eqnarray}
H(c) \bullet f_1(c) &=& (c-b) f_1(c+1) \\
H(c) \bullet f_2(c) &=& (c-1) f_2(c+1) 
\end{eqnarray}
and
\begin{eqnarray}
B(c+1) \bullet f_1(c+1) &=& c^2 f_1(c) \\
B(c+1) \bullet f_2(c+1) &=& c^2\frac{b-c}{1-c} f_2(c).
\end{eqnarray}
The operators $H(c)$ and $B(c+1)$ give 
a left $D$-isomorphism among
$D/I(\beta)$ and $D/I(\beta+e_2)$.
\end{example}

We are interested in the following problem
to apply for our comrehensive restriction algorithm;\\
{\bf Problem}: Find up-step and down-step
contiguity operators that give isomorphisms among left $D$-modules
under a restriction of parameter space or
among $M(\beta)$ and $M(\beta+h)$ where $h$ is any shift direction.

Suppose we reparametrize $\beta$ as
$\beta_1=L_1(\beta'_1, \ldots, \beta'_m),
 \ldots,
 \beta_d=L_d(\beta'_1, \ldots, \beta'_m)
$
where $L_i$ are a linear forms of $\beta'$.
We call this reparametrization 
{\it a restriction of parameter space}\/.
We regard $\beta'$ as a new $\beta$.
For example,
$\beta_1=\cdots=\beta_d = \beta'_1$
is a restriction of parameter space
and $\beta'_1$ is regarded as a new $\beta$.
Our problem is to find an up-step operator
and a down-step operator, which give an isomorphism,
with respect to $\beta'_1$.

How do we find these up-step and down-step contiguity operators?
There are several methods to find contiguity operators
for hypergeometric systems.
In the field of special functions, it is customary to refer to left isomorphisms between rational Weyl algebras $R_n$ (the ring of linear partial differential operators of $n$-variables with rational function coefficients) 
as ``contiguity relations''. 
It should be noted that existing research on contiguity relations,
such as that mentioned below, also encompasses cases of this type.
Here are a list of them.
\begin{enumerate}
\item For given an up-step or a down-step operator, deriving
an down-step operator or an up-step  operator respectively
by Gr\"obner basis  \cite{takayama-1989}, \cite{SST-ip}, \cite{nakayama-2024}.
\item Finding contiguity operators by utilizing the middle convolution
and some other operators for rigid systems \cite{oshima-2012}.
\item Finding $D_n$-isomorphism among GKZ hypergeometric systems \cite{saito-isom}.
\item Finding $R_n$-isomorphism by finding rational solutions of a system of
linear differential equations.
\item Finding $D_n$-isomorphism of Horn type hypergeometric systems by restricting isomorphisms of GKZ hypergeometric systems.
\end{enumerate}
Each method has advantages and disadvantages.
We briefly explain first three known methods by examples.
For general description of these method, please refer to the cited papers above.
Last two methods are new and we will give general 
descriptions in next sections together with examples.

\subsection{Deriving down(up)-step operator for a given up(down)-step contiguity operator }
\label{sec:find_contiguity_by_sygyzy}

Suppose that we are given an up-step (resp. a down-step) operator $H$.
The down-step (resp. the up-step) operator can be constructed by a Gr\"obner basis computation in the ring of differential operators 
when parameters are generic numbers \cite{takayama-1989}, \cite{nakayama-2024}.
Let us explain this method by an example.

\begin{example}\rm \label{ex:gauss_contiguity}
The Gauss hypergeometric equation in terms of Euler operator is
\begin{equation}  \label{eq:gauss_by_Euler_op}
L(\aa,\bb,\cc) \bullet f=0, \quad
 L(\aa,\bb,\cc) = \theta_x(\theta_x+\cc-1) - x (\theta_x+\aa) (\theta_x+\bb).
\end{equation}
Put $H_\aa(\aa)=\theta_x+\aa$.
By the relation 
$x (\theta_x+\aa+1) = (\theta_x+\aa) x $ in $D$,
we have 
\begin{equation} \label{eq:trivial_contiguity}
 L(\aa+1,\bb,\cc) H_\aa(\aa) = H_\aa(\aa) L(\aa,\bb,\cc) \in D L(\aa,\bb,\cc).
\end{equation}
Therefore, the operator $H_\aa(\aa)$
is an up-step operator with respect to $\aa$.
Suppose $H_\aa(\aa)$ gives an isomorphism among $D/L(\aa,\bb,\cc)$
and $D/L(\aa+1,\bb,\cc)$.
Since the inverse of $H_\aa(\aa)$ is a down-step operator 
$B_\aa(\aa+1)$, 
the relation
\begin{equation} \label{eq:BH_relation}
B_\aa(\aa+1) H_\aa(\aa) - 1 \equiv 0 \quad {\rm mod}\  D L(\aa,\bb,\cc)
\end{equation}
holds.
In other words, the down-step operator $B_\aa(\aa+1)$ 
can be obtained by solving the inhomogeneous 
syzygy equation in $D$
\begin{equation} \label{eq:syzygy_gauss}
 -1 + s_1 H_\aa(\aa) + s_2 L(\aa,\bb,\cc)=0
\end{equation}
where $s_1, s_2$ are unknown elements in $D$
and $s_1$ is $B_\aa(\aa+1)$.
There are several algorithms solving inhomogeneous
syzygy equations.
In this case, computing the Gr\"obner basis
of $(H_\aa(\aa),1), (L(\aa,\bb,\cc),0)$ in $D^2$
by the POT order solves the syzygy equation \cite{nakayama-2024}.
The Gr\"obner basis contains an element
$$
(\aa(\cc-\aa-1), x(1-x) \pd{x}-\bb x-\aa+\cc-1)
= c_1 (H_\aa(\aa),1) + c_2 (L,0), \quad c_i \in D,
$$
which implies that
$
B_\aa(\aa+1)=\frac{1}{\aa(\cc-\aa-1)}
 x(1-x) \pd{x}-\bb x-\aa+\cc-1
$.
\end{example}

Since classical hypergeometric systems have 
either a trivial up-step operator or 
a down-step operator
as in (\ref{eq:trivial_contiguity}),
we can obtain any contiguity operator for an integral
shift for generic values of parameters
by a composition and the method of this section. 

We call
up-step operators $H_i(\beta)$ and 
down-step operators $B_i(\beta)$ 
{\it atomic contiguity operators}\/.
When they give isomorphisms, a composite of them
also gives an isomorphism.
However, a restriction in the parameter space 
of the composite 
does not always give an isomorphism.

\begin{example}\rm \label{ex:rest_of_isom_does_not_give_isom}
We denote $x_1$ by $x$ and $D_1$ by $D$.
We consider the Gauss hypergeometric system
$D/DL(\aa,\bb,\cc)$.
The following operators are atomic contiguity operators.
\begin{eqnarray}
H_a(\aa,\bb,\cc)&=&  {x}\pd{x} +  {\aa }, \\
B_a(\aa,\bb,\cc)&=&    -  {x}  (  {x}- 1)\pd{x}
  -   (    {\bb }  {x}+ {\aa }- {\cc }), \\
H_c(\aa,\bb,\cc)&=&  {x}  (  {x}- 1)\pd{x}^{2}
 +  (   (    {\aa }+ {\bb }- {\cc }+ 2)  {x}- 1)\pd{x}
 +  (     (   {\bb }- {\cc }+ 1)  {\aa }+  (  - {\cc }+ 1)  {\bb }+  {\cc }^{ 2} - {\cc }), \\
B_c(\aa,\bb,\cc)&=& {x}\pd{x} +  {\cc }- 1
\end{eqnarray}

We compose them as 
$$ 
   H:=H_\cc(\aa+1,\bb,\cc) H_\aa(\aa,\bb,\cc), \    
   B:=B_\cc(\aa,\bb,\cc+1) B_\aa(\aa+1,\bb,\cc+1). 
$$
Reducing $B$ by $DL(\aa+1,\bb,\cc+1)$, we obtain
\begin{equation} \label{eq:isom-B-gauss1}
{\bar B}=
(  {\aa }- {\cc })   ({x}  (  {x}- 1)\pd{x}
 +   {\bb }  {x}- {\cc }). 
\end{equation}
It gives an isomorphism among
$D/DL(\aa+1,\bb,\cc+1)$ and
$D/DL(\aa,\bb,\cc)$ for generic values of
parameters.
Note that ${\bar B}$ can be divided by $\aa-\cc$. 
When we restrict $B$ to $\aa=\cc$,
we have
$B'=    -  {x}^{ 2}   (  {x}- 1)\pd{x}^{2}
  -    {x}  (    (   {\bb }+ {\cc }+ 2)  {x}- {\cc }- 1)\pd{x}
  -   (  {\cc }+ 1)  ( {\bb }) {x}
$,
which belongs to the left ideal $DL(\cc+1,\bb,\cc+1)$.
This means that $B'$ does not give an isomorphism.
On the other hand, we can see
$D/DL(\cc+1,\bb,\cc+1)$
and
$D/DL(\cc,\bb,\cc)$
are isomorphic
for generic complex numbers $\bb, \cc$
by 
\begin{equation} \label{eq:gauss-h-isom1}
{\bar H}=
  \frac{1}{b(\bb,\cc)} \left( (  {\cc }- 1)  (  {x}- 1)\pd{x}
  -   {\cc }  (  {\cc }- 1) \right)
\end{equation}
and
by ${\bar B}/(\aa-\cc)$
where $b(\bb,\cc)=\cc(\cc-1)(\cc-\bb)$.

This observation shows that a restriction 
of a composite of atomic contiguity operators, 
which gives an isomorphism
for generic values of parameters,
does not always give an isomorphism.
However, dividing a factor like $\aa-\cc$ might
give an isomorphism as we have seen above.
\end{example}


\subsection{Finding contiguity operator for rigid systems} \label{sec:oshima-rigid}
Let us briefly explain a method to construct contiguity relations
given in \cite[Sec 3.2, Chap 11]{oshima-2012} by an example.
We will construct a contiguity relation with respect to $c$ 
for the hypergeometric function
$$
   f(a,b,c;z)=\frac{1}{\Gamma(a+1)} \int_0^1 (z-x)^a x^b (1-x)^c dx.
$$
It satisfies the Gauss hypergeometric equation
(\ref{eq:gauss_by_Euler_op})
$L(-a,-a-b-c-1,-a-b) \bullet f = 0$.
Put
$\phi = x^b (1-x)^c $.
Put $\phi_+=(1-x) \phi$.
Applying $\pd{x}$ to it, we have
$$ \pd{x} \bullet \phi_+ = (\pd{x} -(x \pd{x}+1))\bullet \phi. $$
Now, we apply the fractional derivative
$\pd{x}^{-\mu}$, $\mu=a+1$ to the both sides.
Note that we have the formula
$$ \pd{x}^{-\mu} x \pd{x} = (x \pd{x}-\mu) \pd{x}^{-\mu} $$
or 
$$ \pd{x}^{-\mu} x \pd{x} = \pd{x}^{-\mu} x \pd{x} \pd{x}^\mu \pd{x}^{-\mu} 
  = {\rm Ad}(\pd{x}^{-\mu}) (x \pd{x}) \pd{x}^{-\mu}, \quad
   {\rm Ad}(f)L := f^{-1} L f
$$
in the ring of factional differential operators\footnote{We have no rigorous definition of this ring. The term is used as an intuitive wording.} \cite[Secs 1.2, 1.3]{oshima-2012}.
Moreover, we have
$$
   \pd{x}^{-\mu} \bullet \varphi(x) := I^\mu(\varphi):=\frac{1}{\Gamma(\mu)}\int_c^x \varphi(t) (t-x)^{\mu-1} dt
$$
where $c$ is suitably chosen.
This action gives a left module structure to the ring of fractional 
differential operators and a space of holomorphic functions.
By utilizing these relations, we have
\begin{eqnarray}
 \pd{x}^{-\mu} \pd{x} \bullet \phi_+ &=& \pd{x}^{-\mu} (\pd{x} -(x \pd{x}+1))\bullet \phi \\
 \pd{x} \pd{x}^{-\mu} \bullet \phi_+ &=& (\pd{x} -(x \pd{x}-\mu+1))\pd{x}^{-\mu}  \bullet \phi \\
 \pd{x} \bullet I^{\mu}(\phi_+) &=& (\pd{x} -(x \pd{x}-\mu+1)) \bullet I^\mu(\phi) 
\end{eqnarray}
Thus, changing the variable $x$ by $z$, we have
\begin{equation}  \label{eq:conti0}
  \pd{z} \bullet f(a,b,c+1;z) = (\pd{z}-(z \pd{z} - a)) \bullet f(a,b,c;z)
\end{equation}
The function $f(a,b,c;z)$ satisfies the ODE
\begin{equation}
  L \bullet f(a,b,c;z)=0, \quad
  L=\theta_z (\theta_z-a-b-1)-z (\theta_z-a)(\theta_z-a-b-c-1)
\end{equation}
where $\theta_z = z \pd{z}$.

There exist differential operators $r_3, r_4$ such that
$r_3 \pd{} - 1 = r_4 L(a,b,c+1)=0$,
because $L$ is irreducible for generic values of $a,b,c$.
In fact, 
$r_3 = ((z-z^2)\pd{z}+(2a+b+c+1) z-a-b)/(a (a+b+c+2))$
and 
$r_4=L/(z a (a+b+c+2))$.
Applying $r_3$ to (\ref{eq:conti0}), we have
$$
  f(a,b,c+1;z)= r_3 (\pd{}-(z \pd{z} - a)) \bullet f(a,b,c;z).
$$
Reducing $r_3 (\pd{z}-(z \pd{z}-a))$ by $L(a,b,c)$,
we obtain 
\begin{equation} \label{eq:contiguity_for_integral}
\frac{z (1-z) \pd{z}+a z+(c+1)}{a+b+c+2} \bullet f(a,b,c;z)=f(a,b,c+1;z)
\end{equation}
which is a contiguity relation.

Note that when $b=-1$, the operator $L$ is factored as
\begin{equation}
(\theta_z - z (\theta_z-a-c))(\theta_z-a)
\end{equation}
and then it is not irreducible.
However, we are lucky for the case $b=-1$ 
that the inverse $r_3$ of $\pd{z}$ exists
and the method above works for this degenerate case.
When $a=0, b=-1$, there is no inverse of $\pd{z}$ modulo $L$,
because the left ideal generated by $\pd{z}$ and $L$ is the principal
ideal generated by $\pd{z}$.
The method of \cite{oshima-2012} does not give a contiguity relation
for this case. 
Note that a different approach gives the isomorphism.  See Example \ref{ex:by-rational-sol}, (\ref{eq:A-case-a0-bm1}), and Section \ref{sec:restriction_of_A_hg}.
The contiguity derived by methods above agrees with (\ref{eq:contiguity_for_integral}) restricted to $a=0$ and $b=-1$. 
The agreement seems to be a coincidence.
As we have seen in Example \ref{ex:rest_of_isom_does_not_give_isom}
it is not always possible to obtain an up-step or a down-step operator by a restriction of parameters.


Finally, we note two things.

It follows from the relation (\ref{eq:contiguity_for_integral}) and the comparison of the constant term 
that the contiguity relation for hypergeometric series
$g(a,b,c;z):=F(-a,-a-b-c-1,-a-b;z)$ is
$$
\frac{z (1-z) \pd{z}+a z+(c+1)}{c+1} \bullet g(a,b,c;z)=g(a,b,c+1;z).
$$

The Riemann scheme of the ODE $L f = 0$ is
\begin{equation}
 \left\{
  \begin{array}{ccc}
   x=0 & x=1 & x=\infty \\
    0  &   0 & -a \\
  a+b+1&a+c+1&-a-b-c-1 \\
  \end{array}
 \right\}
\end{equation}

\subsection{Finding isomorphisms among GKZ hypergeometric systems}

Let $A=(a_{ij})$ be a $d \times n$ matrix with integer entries.
We assume the column vectors $a_i$ of $A$ generates the lattice $\ZZ^d$.
Let $\beta \in \CC^d$ be a parameter vector.
The the {\it GKZ left ideal} 
$H_A(\beta)$ 
is generated by
\begin{eqnarray}
&& \sum_{j=1}^n a_{ij} x_j \pd{j} - \beta_i, \quad i=1, \ldots, d \label{eq:gkz_euler} \\
&& \pd{}^u - \pd{}^v, \quad Au=Av,\  u, v \in \NN_0^n.  \label{eq:I_A}
\end{eqnarray}
We denote by $M_A(\beta)$ the left $D$-module $D_n/H_A(\beta)$.
The ideal  generated by (\ref{eq:I_A})
in the polynomial ring $\CC[\pd{1}, \ldots, \pd{n}]$
is called the {\it affine toric ideal} and denoted by $I_A$.
Depending on the context, 
the {\it GKZ (hypergeometric) system} refers either to the system of equations $H_A(\beta)$ 
or to the $D$-module $M_A(\beta)$.

Mutsumi Saito \cite{saito-isom} 
gave a necessary and sufficient condition that  
given GKZ hypergeometric systems\footnote{Although the term ``GKZ hypergeometric system'' appears in many recent publications, the original name was ``A-hypergeometric system''.}
are isomorphic.
He also gave a method to construct an isomorphism among
isomorphic GKZ hypergeometric systems.

Let us see his construction with an example.
Consider  
\begin{equation}  \label{eq:A-of-gauss1}
 A=\left(\begin{array}{cccc}
 1 & 0 & 0 &-1 \\
 0 & 1 & 0 & 1 \\
 0 & 0 & 1 & 1 \\
\end{array} \right).
\end{equation}
and a parameter shift 
\begin{equation}
\chi=\chi_+ - \chi_-, \ \chi_+=(1,0,0)^T, \chi_-=(0,1,0).
\end{equation}
We have
$\chi_+ = A u$, $u=(1,0,0,0)^T$ and
$\chi_- = A v$, $v=(0,1,0,0)^T$.
The monomial ideal $M_\chi=\langle \pd{}^u \,|\, Au \in \chi + \NN_0 A \rangle$
 \cite[(4.13)]{saito-isom} is generated 
by $\pd{1}, \pd{3}$.
We will give an algorithm to construct the monomial ideal $M_\chi$
in subsection \ref{sec:b-chi}.
The $b$-ideal $B_\chi$ for the direction $\chi$  is
\begin{equation}
\left( \left\langle \sum_{j=1}^n a_{ij} \theta_j - s_i \,|, i=1, \ldots,d \right\rangle
  + {\tilde M}_\chi \right) \cap \QQ[s_1, \ldots, s_d].
\end{equation}
Here ${\tilde M}_\chi$ is the distraction \cite{SST} of $M_\chi$.
Then, we can see that the $b$-ideal $B_\chi$ \cite[(4,14)]{saito-isom}
is generated by $b(s)=s_1+s_3$.
We want to construct an operator $E$ such that
\begin{equation}
  E \pd{}^u = b(\beta) \pd{}^v \ {\rm mod}\, H_A(\beta)
\end{equation}
where we regard $\beta$ as indeterminates.
We may regard $E$ as an inverse operator of $\pd{}^{u-v}$.
Saito \cite[Alg 4.2]{saito-isom} gives an algorithm to construct $E$
from the $b$-ideal $B_\chi$.

The following algorithm \ref{alg:simple_alg_to_find_E_b} to find $b$ and $E$ will be simpler for small examples.
Compute Gr\"obner basis in the free module in $D^2$ of
$(\pd{1},1)$, $\{ (\ell,0)\,|\, \ell \in H_A(\beta) \}$
with the POT order such that $x_1,x_2,x_3,x_4, \pd{1},\pd{3},\pd{4} \succ \pd{2}, \beta_1, \beta_2, \beta_3$ \cite{nakayama-2024}.
The Gr\"obner basis contains an element
$((\beta_1+\beta_3)\pd{2},x_1\pd{2}+x_3\pd{3})$.
Then, we have $E=x_1\pd{2}+x_3\pd{4}$.
Let $f(\beta;x)$ be a solution of $H_A(\beta)$.
Then, we have
$ E \pd{1} \bullet f(\beta;x)= (\beta_1+\beta_3) \pd{2} \bullet f(\beta;x)$.
Since $\pd{i} \bullet f(\beta;x)=f(\beta-a_i;x)$ (modulo non-zero constant factor), we have
\begin{equation}
 E \bullet f(\beta-e_1;x)= (\beta_1+\beta_3) f((\beta-e_1)+e_1-e_2;x).
\end{equation}
In other words, $E$ gives a up-step operator for $\chi=e_1-e_2$.
Note that $E$ gives an isomorphism of corresponding $D$-modules
under some conditions.

By setting $\beta=(c-1,-a,-b)$, solutions of this GKZ hypergeometric
system can be written by the Gauss hypergeometric function
${}_2F_1(a,b,c;z)$.
Assume $a=c$. Then the restriction of $E$ (see Section \ref{sec:restriction_of_A_hg}) gives a contiguity
for the integer shift of $c$ for ${}_2F_1(c,b,c;z)$.

\begin{algorithm}[H]\rm 
\footnotesize
\caption{Finding $E$ and $b$}  \label{alg:simple_alg_to_find_E_b} \ 
\begin{itemize}
\item {\bf Require}: generators $\ell_1, \ldots, \ell_m$ of $H_A(\beta)$. $\pd{}^u$, $\pd{}^v$
where $u, v \in \NN_0^d$ and their supports are disjoint and
$D/H_A(\beta-Au)$ and $D/H_A(\beta-Av)$ are isomorphic. 
\item {\bf Ensure}: $E \in D$ and $b \in \CC[\beta]$ such that
$E \pd{}^u = b \pd{}^v$ modulo $H_A(\beta)$.
\end{itemize}
\begin{enumerate}
\item Compute a Gr\"obner basis $G$ by the POT order of
\begin{equation}
 \left( \begin{array}{c}\pd{}^u\\ 1\\ 0\\ 0\\ \cdot\\ \cdot\\ \cdot\\ 0\end{array} \right),
 \left( \begin{array}{c}\ell_1\\ 0\\ 1\\ 0\\ \cdot\\ \cdot\\ \cdot\\ 0\end{array} \right),
 \left( \begin{array}{c}\ell_2\\ 0\\ 0\\ 1\\ \cdot\\ \cdot\\ \cdot\\ 0\end{array} \right),
  \cdots,
 \left( \begin{array}{c}\ell_m\\ 0\\ 0\\ 0\\ \cdot\\ \cdot\\ \cdot\\ 1\end{array} \right)  \in D^{m+2}
\end{equation}  
The tie breaker $\prec$ of the POT order is
$\beta \prec (\mbox{$\pd{i}$'s in the support of $\pd{}^v$}) \prec 
(\mbox{other variables})$.
\item Find an element of the form
$(b \pd{}^w,c_0, c_1, \ldots,c_m)^T$ such that $\pd{}^w | \pd{}^v$ 
in the Gr\"obner basis $G$.
\item Put $E=\pd{}^{v-w} c_0$ and return $E$ and $b$.
\end{enumerate}
\normalsize
\end{algorithm}
Note that each element of $\beta$ may be degree $1$ or $0$ polynomials
of indeterminates.
For example, $\beta=(-c,-c,1,c-1,c'-1)$ is OK and
$\CC[\beta]$ means $\CC[c,c']$.

The correctness of this algorithm can be shown as follows.
The existence of $E$ and $b$ is proved in \cite{saito-isom}.
Therefore, the Gr\"obner basis of $\pd{}^u$ and $\ell_i$'s
must contain an element of the form $\tilde{b} \pd{}^w$, $\tilde{b} \in \CC[\beta]$ 
whose leading term ${\rm in}_\succ(\tilde{b}) \pd{}^w$
divides the leading term ${\rm in}_\succ(b) \pd{}^v$.
Note that $E$ and $b$ are not unique in general. 
We denote $\tilde{b}$ by $b$ in the sequel.
Since the Gr\"obner basis is computed by the POT order,
we have
$ b \pd{}^w = c_0 \pd{}^u + \sum_{i=1}^m c_i \ell_i$
where $c_i \in D$.
Applying $\pd{}^{v-w}$, we obtain the output.

Let us consider a degenerate case of
$\beta=(\beta_1,0,1)$
for our $A$ (\ref{eq:A-of-gauss1}). 
It stands for the case $a=0, b=-1, c-1=\beta_1$,
which were considered in Section \ref{sec:oshima-rigid}.
By applying the algorithm of \cite{saito-isom}, we have
\begin{equation}
 U \pd{1} = \beta_1(\beta_1+1), \ 
 U= -(x_1 x_4-x_2x_3)\pd{4}+(\beta_1+1)x_1,
\end{equation}
which gives an isomorphism of $D$ modules 
\begin{eqnarray} \label{eq:A-case-a0-bm1}
M_A((\beta_1,0,1)) \ni \ell &\longmapsto& \ell \frac{U}{\beta_1(\beta_1+1)} \in  M_A((\beta_1-1,0,1)), \\ 
M_A((\beta_1-1,0,1)) \ni \ell &\longmapsto& \ell\pd{1} \in  M_A((\beta_1,0,1)) \nonumber
\end{eqnarray}
when $\beta_1(\beta_1+1) \not= 0$.

\section{Finding contiguity operators by finding rational solutions} \label{sec:rat-sol}

We consider a family of holonomic $D$-modules $M(\beta)=D/I(\beta)$
where $I(\beta)$ is a left ideal of $D$.
Let $\ell_1, \ldots, \ell_m$ are generators of $I(\beta+e_i)$.
Then, the condition (\ref{eq:H-relation}) satisfied by the up-step operator $H_i$  
is equivalent to
\begin{equation}  \label{eq:H-relation2}
  \ell_j \, H_i(\beta) \in I(\beta), \ j=1, \ldots, m.
\end{equation}
Let $R_n=\CC(x_1, \ldots, x_n) \langle \pd{1}, \ldots, \pd{n}\rangle$
be the rational Weyl algebra (the ring of differential operators
with rational function coefficients).
Let $\{ s_k \,|\, k=1, \ldots, r\}$ be a set of the standard monomials
with respect to a Gr\"obner basis $G$ of the left ideal $R_n I(\beta)$ in $R_n$.
The set is a basis of 
$R_n/(R_n I(\beta))$ as a vector space over the rational function field
$\CC(x)$ where $\CC(x)$ is an abbreviation of $\CC(x_1, \ldots, x_n)$.
Note that we may assume in $R_n$
\begin{equation}  \label{eq:H-relation2-in-R_n}
  \ell_j \, H_i(\beta) \in R_n I(\beta), \ j=1, \ldots, m.
\end{equation}
Then, the operator $H_i$ can be expressed as
\begin{equation}
 H_i = \sum_{k=1}^r c_k(x) s_k 
\end{equation}
where $c_k(x)$ is an element of $\CC(x)$.
Reducing $\ell_j H_i$ by the  Gr\"obner basis $G$,
we have 
$\sum_{k=1}^r (L_{ij}^k \bullet c_k) s_k$
where $L_{ij}^k \in R_n$.
Then, 
\begin{equation} \label{eq:H-by-rational-sol}
L_{ij}^k \bullet c_k = 0, \ j=1, \ldots, m, \ k=1, \ldots, r
\end{equation}
should hold since $\ell_j H_i$ belongs to $I(\beta)$.
From the above discussion, the problem of finding an up-step operator
$H_i$ has been reduced to the problem of finding a rational solution
$c_k$, $k=1, \ldots, r$
of (\ref{eq:H-by-rational-sol}).
This system can be transformed into an integrable connection
(a Pfaffian system).
An algorithm of finding the rational solutions of an integrable connection
is given in \cite{barkatou-2012}.
We utilize this algorithm to solve (\ref{eq:H-by-rational-sol}).

\begin{example}\rm  \label{ex:by-rational-sol}
Consider the left ideal $I(c)$ generated by
$ \ell = (\pd{x} - (x \pd{x}-c)) x \pd{x} $ in the $D=D_1$ of one variable $x=x_1$.
The set of the standard monomials is $\{1, \pd{x}\}$ and
we set $H=c_0(x) + c_1(x) \pd{x}$.
From (\ref{eq:H-by-rational-sol}), the vector valued function
$F=(c_0',c_0,c_1',c_1)^T$ satisfies the equation
$\frac{dF}{dx}=P F$
where
$$ P=
\left(\begin{array}{cccc}
\frac{   {c}  {x}+ 1}{   {x}^{ 2} - {x}}& 0& 0& 0 \\
 1& 0& 0& 0 \\
 - 2& \frac{ 1}{  {x}- 1}& \frac{   (  - {c}+ 2)  {x}- 1}{   {x}^{ 2} - {x}}& \frac{    (   2  {c}- 2)   {x}^{ 2} +  3  {x}- 1}{    {x}^{ 4} -  2   {x}^{ 3} +  {x}^{ 2} } \\
0& 0&  1& 0 \\
\end{array}\right).
$$
The space of rational functions of this equation is spanned by
$(0,-1,\frac{2x-1}{c+1},\frac{x(x-1)}{c+1})$.
Hence, we have
\begin{equation}
H=\frac{x(1-x)}{c+1}\pd{x}+1.
\end{equation}
By computing the Gr\"obner basis of the left $D$ submodule
generated by $(1,x(1-x)\pd{x}+c+1)$ and $(0,L)$ in $D^2$
with the POT order,
we see that $(-x\pd{x}+c+1,(c+1)^2)$ is in the basis.
Hence, we have
$ (-x \pd{x}+c+1) H - (c+1) \in I(c)$,
which means that
\begin{equation}
B(c+1)=\frac{-x}{c+1}\pd{x} + 1.
\end{equation}

\begin{remark}\rm
For $D_n$-ideals, the set of the standard monomials are not finite.
We can apply this method to a finite subset of the set
with looking for polynomial solutions instead of rational solutions.
If $D_n/H(\beta)$ and $D_n/H(\beta')$ are isomorphic as left $D_n$-modules,
there exists a finite subset $\{ s_k \}$
to express $H_i(\beta)$.
Hence, if two $D_n$-modules are isomorphic,
the modified method above can find contiguity operators
by enlarging the finite subset in finite steps.
\end{remark}

\end{example}

\section{From contiguity operators of GKZ hypergeometric systems to those of Horn type hypergeometric systems} \label{sec:restriction_of_A_hg}

A relation between GKZ hypergeometric systems
and classical hypergeometric systems,
which are called {\it Horn type hypergeometeric systems}, studied categorically in
\cite{berkesch}.
We study a relation of them in terms of restriction of $D$-modules. 

We consider an GKZ hypergeometric ideal $H_A(\beta)$.
We assume the $d \times n$ configuration matrix $A$ is of the form
$\left( E_d\ \vert\  A' \right)$
where $E_d$ is the $d \times d$ identity matrix.
For example, $$
 A=\left(\begin{array}{ccc|c}
 1 & 0 & 0 &-1 \\
 0 & 1 & 0 & 1 \\
 0 & 0 & 1 & 1 \\
\end{array} \right).
$$
satisfies this assumption with $A'=(-1,1,1)^T$.
\begin{theorem}  \label{th:bf_for_A}
Assume $A=\left( E_d \ \vert \ A' \right)$.
\begin{enumerate}
\item The $b$-function (indicial polynomial) along 
$x_1= \cdots = x_d =1$ is $s$.
\item 
The restriction 
\begin{equation}
 D/((x_1-1)D + \cdots (x_d-1)D) \otimes_{D} D/H_A(\beta) 
\end{equation}
is isomorphic to
\begin{equation}  \label{eq:Horn_type_as_restriction}
\frac{D_{n-d}}{D_{n-d} \cap (H_A(\beta)+(x_1-1)D + \cdots (x_d-1)D)} 
\end{equation}
where $D_{n-d}=\CC\langle x_{d+1}, \ldots, x_n, \pd{d+1}, \ldots, \pd{n}\rangle$.
\end{enumerate}
\end{theorem}

\begin{proof}
(1) 
It follows from the assumption of the form of $A$,
the left ideal $H_A(\beta)$ contains the operator
$\theta_i + \sum_{j > d}  a_{ij} \theta_j - \beta_i$.
We change the variables $x_i \rightarrow x_i+1$, $i=1, \ldots, d$,
then this operator becomes
$$\theta_i + \pd{i} + \sum_{j > d}  a_{ij} \theta_j - \beta_i.
$$
The initial term of this operator
with respect to the weight vector
$(-w,w)=(-{\bf 1}_d,{\bf 0}_{n-d}, {\bf 1}_d, {\bf 0}_{n-d})$
is $\pd{i}$
where ${\bf 1}_d$ is a row vector of $d$ ones and
${\bf 0}_{n-d}$ is a row vector of $n-d$ zeros.
Then, the initial ideal with respect to $(-w,w)$ of $H_A(\beta)$ 
with the new coordinates contains
$\pd{1}, \ldots, \pd{d}$.
Therefore
$$\CC[\theta_1+\cdots +\theta_d] \cap {\rm in}_{(-w,w)} (H_A(\beta))$$
contains $s=\theta_1+\cdots +\theta_d$.
Since $H_A(\beta)$ is regular holonomic,
it is specializable and the $b$-function is not constant.
Thus, we have $b(s)=s$.

The statement (2) follows from the restriction algorithm \cite{oaku-1997} 
and the statement (1).
\end{proof}

We call (\ref{eq:Horn_type_as_restriction}) the Horn type hypergeometric system
associated to the GKZ system $H_A(\beta)$.
The $D$-module (\ref{eq:Horn_type_as_restriction})
and the left ideal in the denominator of it
are also denoted by $M_A(\beta)$ and $H_A(\beta)$ respectively
if no confusion arises.

The Gauss hypergeometric system 
is the left $D_1$ module defined by $D_1/(D_1 L)$,
\begin{equation}
L= x(1-x)\p^2 + (\cc - (\aa + \bb + 1)x)\p - \aa \bb
\end{equation}
where $x_1$ is denoted by $x$.
The Appell $F_1$ system is the left $D_2$  module defined by 
$D_2/I_{F_1}$
where $I_{F_1}$ is the left ideal generated by 
(\ref{eq:f1a}), (\ref{eq:f1b}), (\ref{eq:f1c}) 
where $(x_1,x_2)$ is denoted by $(x,y)$.

The Appell $F_2$ system is the left $D_2$ module defined by
$D_2/I_{F_2}(\aa,\bb,\bb',\cc,\cc')$
where $I_{F_2}(\aa,\bb,\bb',\cc,\cc')$ is the left ideal generated by
(\ref{eq:f2a}), (\ref{eq:f2b}).

\begin{theorem} \label{th:restriction_to_gauss_F1_F2}
For any parameter value,
the restriction of the following GKZ hypergeometric systems
defined by $A$ and $\beta$ as
(\ref{eq:gauss-A}), (\ref{eq:appellF1-A}), (\ref{eq:appellF2-A})
to $x_1= \cdots = x_d=1$ are the Gauss hypergeometric system,
the Appell $F_1$ system, and the Appell $F_2$ system respectively
by changing the variable names appropriately, e.g.,
$(x_6,x_7)$ is $(x,y)$ in the case of $F_2$.
\begin{eqnarray}
 A&=&\left(\begin{array}{cccc}
 1 & 0 & 0 &-1 \\
 0 & 1 & 0 & 1 \\
 0 & 0 & 1 & 1 \\
\end{array} \right), \quad
 \beta=(\gamma-1,-\alpha,-\beta)^T  \label{eq:gauss-A} \\
A&=&\left(\begin{array}{cccccc}
1& 0& 0& 0& 1& 1\\
0& 1& 0& 0& 1& 0\\
0& 0& 1& 0& 0& 1\\
0& 0& 0& 1& -1& -1\\
\end{array}\right), \ 
\beta=(-\aa,-\bb,-\bb',\cc-1)^T \label{eq:appellF1-A} \\
A&=&\left(\begin{array}{ccccccc}
1& 0& 0& 0& 0& 1& 1\\
0& 1& 0& 0& 0& 1& 0\\
0& 0& 1& 0& 0& 0& 1\\
0& 0& 0& 1& 0& -1& 0\\
0& 0& 0& 0& 1& 0& -1
\end{array}\right), \ 
\beta=(-\aa,-\bb,-\bb',\cc-1,\cc'-1)^T \label{eq:appellF2-A}
\end{eqnarray}
\end{theorem}

\begin{proof}
Change variables $x_i$ to $x_i+1$ for $i=1, \ldots, d$.
Compute a $(-w,w)$ Gr\"obner basis $G$ for the restriction of 
each GKZ hypergeometric system with the tie breaking block order satisfying
$x_1, \ldots, x_n, \pd{1}, \ldots, \pd{n} \succ \beta_1, \ldots, \beta_d$
(parameters are last).
Computation by a computer program shows that $(-w,w)$ order of each element of $G$
is positive, $0$, or $-1$.
Note that the b-function for the restriction is $s$ by Theorem \ref{th:bf_for_A}.
Then, the restriction is generated by
$ g_{x_1=\cdots=x_d=0}, \{ g \in G\,|\, {\rm ord}_{(-w,w)}(g)=0\}$
and
$ (\pd{i} g)_{x_1=\cdots=x_d=0}, \{ g \in G\,|\, {\rm ord}_{(-w,w)}(g)=-1\},
i=1, \ldots, d$.
Computation by a computer program shows that the restriction agrees with
the corresponding classical hypergeometric systems.
\end{proof}

\subsection{Restriction of a left $D$-homomorphism}

Let $D=D_n$ be the Weyl algebra in $n$ variables and $I$ a left holonomic ideal in $D$.
$b$-function along $x_n=0$ is the monic generator $b(\theta_n)$ of
the principal ideal
$ {\rm in}_{(-w,w)} (I) \cap \CC[\theta_n]$
where $w=(0,\ldots, 0,1)$ and $\theta_n=x_n \pd{n}$.
Assume $k_0$ be the maximal non-negative root of $b(s)=0$.
Let 
\begin{equation}
 F_{k_0} = \sum_{k=0}^{k_0} D_{n-1} \pd{n}^k 
\end{equation}
Then, the restriction algorithm \cite{oaku-1997} gives a Gr\"obner basis $G \subset F_{k_0}$
such that 
$D/(I+x_n D)$ is isomorphic to 
$F_{k_0}/D_{n-1}G$ as the left $D_{n-1}$ module.

The $b$-function plays a cruicial role in the restriction algorithm.
It follows from the definition of the $b$ function that 
there exists an operator $r$ such that
$$
  b(\theta_n)- r \in I, {\rm ord}_{(-w,w)}(r) \leq -1
$$ 
The key identity is 
\begin{equation}  \label{eq:b-reduction}
 \pd{n}^j b(\theta_n) = b(j)\pd{n}^j + (b(\theta_n+j)-b(j))\pd{n}^j + \pd{n}^j r \quad {\rm mod}\, I.
\end{equation}
Note that 
$
(b(\theta_n+j)-b(j))\pd{n}^j \in x_n D
$
and
$ {\rm ord}_{(-w,w)} (\pd{n}^j r) \leq j-1
$.

We define the normal form of $f \in D$ in $D/(I+x_n D)$ as follows.
\begin{enumerate}
\item Remove all $\pd{n}^j$ ($j > k_0$) and $x_n$ in $f$ 
by (\ref{eq:b-reduction}) modulo $I+x_nD$.
 The result ${\tilde f}$ is in $F_{k_0}$. 
\item Compute the normal form ${\tilde f}$ by the Gr\"obner basis $G$.
 We denote the result by ${\bar f}$.
\end{enumerate}

Assume $\ell \in D$ defines a left $D_n$-morphism among $D/I$ and $D/I'$
by $D/I \ni [f] \mapsto [f \ell] \in D/I'$.
Since it is well-defined, we have
$I \ell \subset I'$.
This morphism induces the left $D_{n-1}$-morphism
\begin{equation}  \label{eq:restriction-morphism}
 D/(I+x_n D) \ni [f] \mapsto [f \ell] \in D/(I'+x_nD)
\end{equation}
It is well-defined because
$(I+x_nD) \ell \subset I'+x_n D$.
The maximal integral root of the $b$-function of $I'$ along $x_n=0$
is denoted by $k'_0$ and
the Gr\"obner basis
obtained by applying the restriction algorithm to $I'$
by $G'$.

\begin{proposition}
Assume $k_0=k'_0=0$.
Then, the morphism (\ref{eq:restriction-morphism}) is given by
\begin{equation}
 D_{n-1}/D_{n-1} G \ni [f] \mapsto [f {\bar {\ell}}] \in  D_{n-1}/D_{n-1} G'
\end{equation}
where the normal form ${\bar {\ell}}$ is taken in $I'+x_n D$.
\end{proposition}

{\it Proof\/.}
Since $F_0=D_{n-1}$, $f$ does not contain the variables $x_n$ and $\pd{n}$.
Then we have 
 $f {\bar \ell} - f \ell  \in x_n D + I'$
from ${\bar \ell} - \ell = x_n c_1 + c_2$, $c_1 \in  D$, $c_2 \in I'$.
Note that $f {\bar \ell} \in D_{n-1}$.
\qed

We note that these results can be easily generalized to 
the case of the restriction to $x_m = x_{m+1} = \cdots = x_n = 0$.

\subsection{Restriction of isomorphisms of GKZ hypergeometric systems to those of classical hypergeometric systems --- restriction of Saito's isomorphism}

We can obtain an isomorphisms among a contiguous family of classical
hypergeometric system such as the Gauss hypergeometric system and
Appell hypergeometric system $F_2$ 
by applying the restriction algorithm to isomorphisms constructed by
M.Saito \cite{saito-isom}
as long as the maximal integral root of the $b$-function for the restriction
is $0$.
Note that this method works for any degenerated parameters.

The general algorithm of computing the restriction of a homomorphism
can be described in a simple form for the GKZ system when $A=(E_d,A')$
and the restriction is that to $x_1=\cdots=x_d=1$.

\begin{algorithm}[H]\rm 
\footnotesize
\caption{LR(left-right)-reduction}  \label{alg:LR_reduction}
\begin{algorithmic}[1]
\Require Rules $\pd{i} \rightarrow \ell_i \in D$, $i=1, \ldots, d$. An element $\ell \in D$. 
\Ensure ${\bar \ell}$ such that
${\bar \ell} = \ell$ modulo $I + \sum_{i=1}^d x_i D$
where $I$ is the left ideal in $D$ generated by $\pd{i}-\ell_i$, $i=1, \dots, d$. 
\Repeat
\State $\ell \leftarrow \ell |_{x_1=\cdots=x_d=0}$; 
\State Choose a term of the form $t:=c x^\alpha \pd{}^\beta \pd{i}$, $c \in K$ in $\ell$
and rewrite 
$$\ell \leftarrow c x^\alpha \pd{}^\beta \ell_i + (\ell - t)$$ 
\Until {$(\mbox{there is no term divided by $\pd{i}$, $i=1, \ldots, d$})$}
\State \Return $\ell$ as ${\bar \ell}$.
\end{algorithmic}
\normalsize
\end{algorithm}

It is easy to see ${\bar \ell}$ satisfies the output condition
when the algorithm stops.
For the GKZ system with $A=(E_d,A')$, we firstly make the change of variables
$J=H_A(\beta) |_{x_i \rightarrow x_i+1, i=1, \ldots, d}$
and use the rules
\begin{equation}
\pd{i} \rightarrow \ell_i, \ 
\ell_i= - x_i \pd{i} - \sum_{i=d+1}^n a_{ij} x_j \pd{j} + \beta_i,
\quad i= 1, \ldots, d.
\end{equation}
Note that $\pd{i}-\ell_i$ belongs to the GKZ ideal $J$.
The LR-reduction choosing $t$ by the lexicographic order 
$\pd{1} \succ \pd{2} \succ \cdots$ stops for this case
because $\ell_i$ contains only the term -$x_i \pd{i}$ and
other terms of $\ell_i$ do not contain $\pd{k}$, $k=1, \ldots, d$.
More precisely, it can be proved as follows.
Consider the ${\rm degree}(c x^p \pd{}^q \pd{i},\pd{i})$
be the degree of $c x^p \pd{}^q \pd{i}$ with respect to $\pd{i}$.
The degrees of all terms in $x^p \pd{}^q \ell_i |_{x_1=\cdots=x_d=0}$
are strictly smaller than the original degree.
We use the lexicographic order $\pd{1} \succ \pd{2} \succ \cdots$
to choose the term $t$.
Then the degree of the leading term decreases strictly in a finte steps.
Then, the LR-reduction stops.

\begin{example} \rm
Let us consider the GKZ hypergeometric system $H_A(\beta)$ for
$A=\left(\begin{array}{cccc}
 1&0&0& -1 \\
 0&1&0&  1 \\
 0&0&1&  1 \\
\end{array}\right)$
and $\beta=(c-1,-a,-b)$.
The column vectors of $A$ is denoted by $a_1, a_2, a_3, a_4$.
We denote $D_4/H_A(\beta)$ by $M_A(\beta)$.
We will restrict $H_A(\beta)$ to $x_1=x_2=x_3=1$.
In other words, we consider
\begin{equation}
\frac{D_4}{H_A(\beta)+(x_1-1)D_4+(x_2-1)D_4+(x_3-1)D_4}.
\end{equation}
The $b$-function $b(s)$ along $x_1=x_2=x_3=1$ is
$s$ and then the maximal integral root is $0$ for any value of $\beta$.
When $b(\beta+a_1)=(\beta_1+1+\beta_2)(\beta_1+1+\beta_3)$ is not zero,
$\pd{1}$ gives an isomorphism 
\begin{equation} \label{eq:gauss-isom}
M_A(\beta) \in [f] \mapsto [f \pd{1}] \in M_A(\beta+a_1)
\end{equation}
and the inverse of $\pd{1}$ is 
\begin{equation}
U_1 = x_2 x_3 \pd{4} + x_1 x_3 \pd{3} + x_1 x_2 \pd{2} + x_1^2 \pd{1} + x_1
\end{equation}
divided by $b(\beta+a_1)$.
See \cite{SST-ip} and \cite{saito-isom} as to algorithms.
Let us compute the normal form ${\bar {U_1}}$.
To do this, we make the change of variables
$x_i \rightarrow x_i+1$ ($i=1,2,3$) in $H_A(\beta)$ and 
consider the restriction to $x_i=0$ ($i=1,2,3$).
The operators $U_1$ is 
\begin{equation} \label{eq:new_U_1}
 (x_2+1)( x_3+1) \pd{4} + (x_1+1)( x_3+1) \pd{3} 
+(x_1+1)( x_2+1) \pd{2} + (x_1+1)^2 \pd{1} + x_1+1
\end{equation}
and first order operators in $H_A(\beta)$ is
\begin{eqnarray}
&&(x_1+1)\pd{1} - x_4 \pd{4} - \beta_1, \\
&&(x_2+1)\pd{2} + x_4 \pd{4} - \beta_2, \\
&&(x_3+1)\pd{3} + x_4 \pd{4} - \beta_3. \\
\end{eqnarray}
Then,
\begin{eqnarray}
\pd{1} &\rightarrow& x_4 \pd{4} + \beta_1, \\
\pd{2} &\rightarrow& -x_4 \pd{4} + \beta_2, \\
\pd{3} &\rightarrow& -x_4 \pd{4} + \beta_3 \\
\end{eqnarray}
are reduction rules (\ref{eq:b-reduction}) obtained by the $b$-function.
Applying these rules to (\ref{eq:new_U_1}) and remove elements in 
$x_1D_4+x_2D_4+x_3D_4$,
we obtain
\begin{equation}
{\bar {U_1}} =
 \pd{4} + (-x_4 \pd{4} + \beta_3)
+(-x_4 \pd{4} + \beta_2)  + (x_4 \pd{4} + \beta_1) + 1
\end{equation}
Replacing $\beta_i$'s by $a,b,c$, we have the contiguity operator
\begin{equation}
\frac{1}{(c-a)(c-b)}\left(
 (1-x_4) \pd{4} - a - b + 1 \right).
\end{equation}
Note that $\bar {\pd{1}}$ is $x_4 \pd{4} + c$.
\end{example}

The isomorphism (\ref{eq:gauss-isom})  holds when
$a=0, b=-1$ and $c(c+1) \not=0$ and then
it induces the isomorphism of the restriction.

If two $A$ hypergeometric systems
$D/H_A(\beta)$ and $D/H_A(\beta')$ are
isomorphic, then the restriction 
of them are isomorphic.

\section{Representatives of Isomorphism Classes} \label{sec:isom}

We consider {\it Horn type} hypergeometric systems obtained by restricting the GKZ hypergeometric system associated with $A=(E_d, A')$ to $x_1=\cdots=x_d=1$.
In the following, we also assume that $\QQ_{\geq 0}A$ is a pointed cone.
Unless there is a risk of confusion, we also denote this Horn type system by $H_A(\beta)$.
For simplicity, we assume $\beta \in \ZZ^d$. 
M. Saito showed that the isomorphism classes of GKZ hypergeometric systems can be described by the set $E_\tau(\beta)$ \cite{saito-isom}. 
Based on this result, we provide an algorithm for finding isomorphism classes based on polyhedral geometry and monomial algebra.
Our algorithm outputs the isomorphism classes as well as the explicit left $D$-module isomorphisms between hypergeometric systems belonging to the same class.
Our algorithm works well for the Gauss hypergeometric system and Appell's $F_1, F_2$ systems, and it is proven to terminate in finitely many steps if $A$ is normal.
As long as we have tried, it works well even when $A$ is not normal.
However, if $A$ is not normal, its termination in finitely many steps has not yet been proven. 
Note that it remains an open question whether the isomorphism of two Horn type hypergeometric systems implies the isomorphism of the corresponding GKZ systems.

\subsection{Algorithm}

\begin{algorithm}[H]
\footnotesize
\caption{\textsc{Representatives}}
\label{alg:rep}
\begin{algorithmic}[1]
\Require $E$: a set of $d-1$ dimensional affine spaces, $param$: a list of parameters
\Ensure $(IsomClass, Contiguity, IsomClassClosure)$
\If{$E = \emptyset$}
    \State $E \gets \textsc{StartingWall}(param)$
\EndIf
\State $L \gets \text{ambient space } \mathbf{Q}^d$
\State \Return \textsc{RepresentativesInAffineSubspace}($L, E, param$)
       \Comment{Call Algorithm \ref{alg:rep_sub}}
\end{algorithmic}
\normalsize
\end{algorithm}

\begin{theorem}  \label{th:correctness}
\begin{enumerate}
\item If Algorithm \ref{alg:rep} terminates in finitely many steps, then this algorithm classifies the Horn type hypergeometric systems $H_A(\beta)$ for $\beta \in \ZZ^d$ into a finite number of isomorphism classes.
\item If $\QQ_{\geq 0} A$ is a normal cone and we take the set of primitive supporting hyperplanes of $A$ as $E$, then this algorithm terminates in finitely many steps.
\end{enumerate}
\end{theorem}

\begin{algorithm}[H]
\footnotesize
\caption{\textsc{RepresentativesInAffineSubspace}}
\label{alg:rep_sub}
\begin{algorithmic}[1]
\Require $L$: an affine subspace, $E$: a set of $d-1$ dimensional affine spaces, $param$: parameters
\Ensure $(IsomClass, Contiguity, IsomClassClosure)$
\State $S \gets \textsc{LinearParametrization}(L)$
\State $dim \gets \dim(L)$
\State $Arr \gets \textsc{FacesOfArrangementOnSubspace}(E, L)$
\State $IsomClass \gets \emptyset$, $Contiguity \gets \emptyset$, $IsomClassClosure \gets \emptyset$
\If{$dim < 1$}
    \State \Return $(Arr, (\emptyset, \emptyset, \text{`0-dim'}), Arr)$
\EndIf
\For{\textbf{each} $C \in Arr$ \textbf{with} $\dim(C) = dim$ } 
    \State $C_{int} \gets \textsc{StrictInteriorFaceForLattice}(C)$
    \State \Comment{C is a face. Try to construct isomorphisms in $C_{int}$}
    \If{$C_{int} = \emptyset$} \textbf{continue} \EndIf
    \If{\textsc{IsBounded}(C)}
        \State Classify $C_{int}$ by exhaustive search and append the data to $IsomClass$, $Contiguity$, $IsomClassClosure$.
        \State \textbf{continue}
    \EndIf
    \State $shifts \gets \textsc{FindChamberBasesAndShifts}(C_{int})$
    \State $C_{cont} \gets \emptyset$
    \For{\textbf{each} $h \in shifts$ } 
        \State $start \gets S$
        \State $target \gets S + h$         
        \Comment{$h$ is a Hilbert basis. Try to construct isomorphism to the direction $h$.}
        \State $walls \gets \textsc{Isom}(S, start, target)$
        \For{\textbf{each} $W \in walls$}
           \State \Comment{No isomorphism is found across $W$ from start to target.}
            \If{$\textsc{Intersection}(W, C_{int}) \neq \emptyset$}
                \If{\textbf{not} \textsc{IsAffineSpaceIncludedInArrangement}(W,L,E)}
                    \State $E'$=\textsc{GetNewHyperplanes}(W,L,E)
                     \Comment{Add new wall and restart}
                    \State \Return \textsc{RepresentativesInAffineSubspace}($L, E', param$) \Comment{Call Algorithm \ref{alg:rep}.}
                \EndIf
            \EndIf
        \EndFor
        \State Append $(\textsc{Contiguity}(start, target), \textsc{Contiguity}(target, start), h)$ to $C_{cont}$
    \EndFor
    \State Append $C_{int}$ to $IsomClass$, $C_{cont}$ to $Contiguity$, $C$ to $IsomClassClosure$
\EndFor
\For{\textbf{each} $C \in Arr$ \textbf{with} $\dim(C) = dim - 1$}
    \State $L' \gets \textsc{GetAffineHullEquations}(C)$
    \State $(I', C', \overline{I}') \gets \textsc{RepresentativesInAffineSubspace}(L', E, param)$
    \State Merge $I', C', \overline{I}'$ into $IsomClass, Contiguity, IsomClassClosure$
\EndFor
\State \Return $(IsomClass, Contiguity, IsomClassClosure)$
\end{algorithmic}
\normalsize
\end{algorithm}

Since the details of this algorithm are technical, let us first explain the outline.
We refer to the $\QQ$-vector space $\QQ^d$ to which the parameter space belongs as the ambient space.
Let $E$ be a set of hyperplanes in the ambient space, and $L$ be an affine subspace within the ambient space.
We consider the arrangement ${\cal A}(E)_L$, which is the arrangement defined by $E$ restricted to $L$.
Let $P$ be a linear parameterization of $L$. The \textsc{RepresentativesInAffineSubspace}($L, E, P$) in Algorithm \ref{alg:rep_sub} performs the classification into isomorphism classes for each face $C$ of the arrangement ${\cal A}(E)_L$ based on the following strategy:
\begin{enumerate}
\item Find a polyhedron $C_{int}$ such that $(\text{interior of } C) \cap \ZZ^d = C_{int} \cap \ZZ^d$.
\item Compute the Hilbert basis of $C_{int} \cap \ZZ^d$, and for its element $h$, compute the condition ($b$-ideal) under which $H_A(P)$ and $H_A(P+h)$ are not isomorphic. It is known that this condition can be written as a union of affine subspaces $W$ \cite[Th 4.5]{saito-isom}.
\item If $C_{int} \cap W$ is empty for all $W$ corresponding to any $h$, then $C_{int}$ is an isomorphism class. If not, add new hyperplanes to $E$ and restart the computation from the beginning.
\end{enumerate}

With the exception of the $b$-ideal computation, we utilize well-known algorithms from hyperplane arrangements and algebra for this process.

The topological structures and the faces of the arrangement are extracted using the classical arrangement construction algorithm detailed by Edelsbrunner et al.\ \cite{edelsbrunner-1986}.
For the systematic enumeration of the faces of the hyperplane arrangement in the parameter space, the reverse search algorithm introduced by Avis and Fukuda \cite{avis-fukuda-1992} works well in practice.
We use a simpler one\footnote{\textsc{FacesOfArrangementOnSubspace}} in our implementation.

To compute the Hilbert basis of the polyhedral cone associated with each chamber, we employ the algorithms implemented\footnote{\textsc{FindChamberBasesAndShifts}} 
in Normaliz \cite{bruns-2010}, \cite{normaliz}.
We note that the computation algorithm of the Hilbert bases for the strict interior of the chambers was proposed by Hemmecke \cite{hemmecke-2002}.

To compute the $b$-ideal\footnote{\textsc{Isom}} used to determine non-isomorphism, we provide an algorithm to compute $B_{-\chi,\chi}$, which appears in Theorem 4.9 of \cite{saito-isom}. 
We will explain this computation algorithm in detail in a later subsection \ref{sec:b-chi}.
Isomorphisms (contiguity relations)\footnote{\textsc{Contiguity}} are
derived by methods explained in previous sections.

The above is an outline of the algorithm, and describing it formally including the technical details yields Algorithms \ref{alg:rep} and \ref{alg:rep_sub}.
We implemented this algorithm using SageMath \cite{sagemath} and its subsystem Risa/Asir \cite{risa-asir}. 
The \textsc{Isom} and \textsc{Contiguity} functions are implemented using Risa/Asir, while the other parts are implemented using SageMath. 
Note that this implementation was developed in discussion with Gemini AI.

We add a note regarding the computational complexity of the algorithm proposed in this paper. 
The primary bottleneck in computation time for this method is not the geometric construction of the hyperplane arrangement in the parameter space, but rather the algebraic steps performed on each chamber---namely, the isomorphism check based on deriving the $b$-ideal involving Gr\"obner basis computation, and the calculation of contiguity relations. 
For enumerating the faces of a hyperplane arrangement comprehensively, highly optimized methods such as the reverse search algorithm by Avis and Fukuda \cite{avis-fukuda-1992} are widely known. 
However, as mentioned above, the overall computation time is dominated by the isomorphism checks of $D$-modules. Therefore, in the implementation of this study, we adopted a simpler and more direct face enumeration algorithm. 
For our target applications, such as GKZ and Horn type hypergeometric systems, as long as the dimension of the parameter space is not excessively high, this simple geometric approach does not compromise overall performance and remains highly practical.

\begin{example}\rm  \label{ex:simple}
Let us look at this algorithm using the simplest example.
The GKZ system $H_A(s)$\footnote{While $\beta$ was a numerical vector, $s$ can be viewed as either a numerical vector or a vector of indeterminates. We omit explaining which one it is when it is clear from the context.} defined by
$A=\left(\begin{array}{cccc}
 1&0&1&0 \\
 0&1&0&1 \\
\end{array}\right)
$
is given by:
$$ x_1 \pd{1}+x_3 \pd{3} - s_1,
   x_2 \pd{2}+x_4 \pd{4} - s_2,
   \pd{2}-\pd{4}, \pd{1}-\pd{3}
$$
Restricting this to $x_1=x_2=1$ yields:
$$ 
  (x_3+1)\pd{3}-s_1, (x_4+1) \pd{4}-s_2
$$
Let us denote the left ideal of the Weyl algebra $D=\QQ \langle x_3, x_4, \pd{3}, \pd{4} \rangle$ generated by these two elements as $I(s)$.
We also define $M(s)=D/I(s)$.
Note that the classical solutions for this system are $(x_3+1)^{s_1}(x_4+1)^{s_2}$.
Here, $d=2$.
Let us execute the algorithm \textsc{RepresentativesInAffineSubspace}($[\ ],[s_1=0, s_2=0],[s_1,s_2]$).
Here $[\ ]$ means $L=\QQ^2$.
Since $L$ is the ambient space, we can directly use $(s_1, s_2)$ as the parameterization for $L$.
The arrangement defined by $s_1 s_2=0$ has four 2-dimensional faces, four 1-dimensional faces, and one 0-dimensional face.
Consider the 2-dimensional face $C=\{s \,|\, s_1 \geq 0, s_2 \geq 0 \}$.
In this case, $C_{int} = \{s \,|\, s_1 \geq 1, s_2 \geq 1 \}$.
The Hilbert basis for $C_{int} \cap \ZZ^d$ is $\{ (1,0), (0,1) \}$.
Let us consider $h=(1,0)$.
We aim to construct an isomorphism between $M(s_1,s_2)$ and $M(s_1+1,s_2)$.
Consider the following morphisms of left $D$-modules\footnote{If we consider the solution space $Hom_D(M(s),{\cal O})$, the direction of the arrows is reversed. This corresponds to contiguity relations in the theory of special functions.}:
\begin{eqnarray*}
M(s_1,s_2) \ni \ell &\mapsto& \ell \frac{1}{s_1+1} \pd{3} \in M(s_1+1,s_2)\\
M(s_1,s_2) \ni \ell (x_3+1) &\leftarrow \! \shortmid& \ell  \in M(s_1+1,s_2)\\
\end{eqnarray*}
Considering the composition
$$ 
M(s_1,s_2) \ni \ell \mapsto \ell \frac{1}{s_1+1} \pd{3} (x_3+1) \in M(s_1,s_2)
$$
we have
$$
\frac{1}{s_1+1} \pd{3} (x_3+1)
= \frac{1}{s_1+1} (x_3+1) \pd{3} + \frac{1}{s_1+1}
\equiv \frac{s_1}{s_1+1} + \frac{1}{s_1+1} = 1 \, {\rm mod}\, I(s_1,s_2)
$$
Even if we reverse the order of composition, we still obtain $1$; hence, this provides an isomorphism if $s_1+1 \not = 0$.
That means the only $Wall$ in this case is $s_1+1 = 0$.
This does not intersect with $C_{int}$.
Considering the $Wall$ when $h=(0,1)$ in a similar manner, we find it is uniquely $s_2+1=0$.
This also does not intersect with $C_{int}$.
Therefore, $C_{int}$ will be registered in $IsomClass$.
The same logic applies to the other 2-dimensional faces.

Next, let us examine the 1-dimensional face $C=\{ s \,|\, s_1 \geq 0, s_2 = 0 \}$.
Since the affine hull of $C$ is $s_2=0$, 
\textsc{RepresentativesInAffineSubspace}($[s_2=0], [s_1=0, s_2=0],[s_1,s_2]$) is executed.
The parameterization $S$ is $(s_1,0)$. 
The arrangement in this case has two 1-dimensional faces and one 0-dimensional face.
Consider the 1-dimensional face $C=\{ s \,|\,  s_1 \geq 0, s_2=0 \}$ of this arrangement.
We have $C_{int} = \{ s \,|\,  s_1 \geq 1, s_2=0 \}$.
The Hilbert basis is $h=(1,0)$.
Consider the following morphisms of left $D$-modules:
\begin{eqnarray*}
M(s_1,0) \ni \ell &\mapsto& \ell \frac{1}{s_1+1} \pd{3} \in M(s_1+1,0)\\
M(s_1,0) \ni \ell (x_3+1) &\leftarrow \! \shortmid& \ell  \in M(s_1+1,0)\\
\end{eqnarray*}
Reasoning identically to before, the walls are $s_1+1=0, s_2=0$.
Because these walls and $C_{int}$ do not intersect, this $C_{int}$ is also registered in $IsomClass$.
Figure \ref{fig:simple_example} illustrates the sets $C_{int}$ output by this algorithm.

Notice that, for example, all $M(s_1,s_2)$ for $s_1, s_2 \in \ZZ_{\geq 0}$ are isomorphic. 
Therefore, it becomes necessary to connect these isomorphic $C_{int}$ regions to each other using the \textsc{MergeIsomClasses} algorithm, which will be described later in subsection \ref{sec:merging}.
\end{example}

\begin{figure}[htbp]
    \centering
    \begin{tikzpicture}[scale=1.2]
        
        \fill[gray!30] (1, 1) rectangle (3, 3);
        \fill[gray!30] (-3, 1) rectangle (-1, 3);
        \fill[gray!30] (1, -3) rectangle (3, -1);
        \fill[gray!30] (-3, -3) rectangle (-1, -1);

        \draw[black, line width=2.5pt] (1, 0) -- (3, 0);
        \draw[black, line width=2.5pt] (-3, 0) -- (-1, 0);
        \draw[black, line width=2.5pt] (0, 1) -- (0, 3);
        \draw[black, line width=2.5pt] (0, -3) -- (0, -1);

        \draw[->, >=stealth, thin] (-3.5, 0) -- (3.5, 0) node[right] {$s_1$};
        \draw[->, >=stealth, thin] (0, -3.5) -- (0, 3.5) node[above] {$s_2$};

        \foreach \x in {-3, -2, -1, 1, 2, 3} {
            \draw (\x, 2pt) -- (\x, -2pt) node[below] {\small $\x$};
        }
        \foreach \y in {-3, -2, -1, 1, 2, 3} {
            \draw (2pt, \y) -- (-2pt, \y) node[left] {\small $\y$};
        }

        \foreach \x in {-3,...,3} {
            \foreach \y in {-3,...,3} {
                \fill (\x, \y) circle (1.2pt);
            }
        }

        \draw[black, very thick, fill=white] (0, 0) circle (2.5pt);

    \end{tikzpicture}
    \caption{An example of isomorphism classes separated by arrangements. The parameter space $(s_1, s_2)$ is partitioned into 2D chambers (shaded), 1D segments (thick lines), and a 0D point at the origin (circle).}
    \label{fig:simple_example}
\end{figure}

Let us explain the data structure used to store hyperplanes and faces in our implementation.
The three basic patterns are as follows:
\begin{enumerate}
\item A hyperplane $c_0 + c_1 s_1 + \cdots + c_d s_d=0$ is represented by $[c_0, c_1, \ldots, c_d]$.
\item An inequality region $c_0 + c_1 s_1 + \cdots + c_d s_d \geq 0$ is also represented by $[c_0, c_1, \ldots, c_d]$.
When this representation is included in a face structure, it denotes an inequality.
An equality $c_0 + c_1 s_1 + \cdots + c_d s_d=0$ within a face structure is represented by the pair $[c_0, c_1, \ldots, c_d]$ and $[-c_0, -c_1, \ldots, -c_d]$.
\item An affine subspace is represented as a list of hyperplane representations.
\end{enumerate}

\begin{example}\rm \label{ex:simple_data_structure}
Taking Example \ref{ex:simple} (the arrangement $s_1s_2=0$), it is represented as follows:
{\footnotesize
\begin{verbatim}
sage: Eq=[[0,1,0],[0,0,1]] # arrangement  
sage: Arr=list_all_arrangement_faces(Eq) 
[[2, [[0, 1, 0], [0, 0, 1]]],  # representation of a face. The leading '2' is the dimension
 [2, [[0, 1, 0], [0, 0, -1]]],
 [2, [[0, -1, 0], [0, 0, 1]]],
 [2, [[0, -1, 0], [0, 0, -1]]],
 [1, [[0, 1, 0], [0, 0, 1], [0, 0, -1]]],
 [1, [[0, -1, 0], [0, 0, 1], [0, 0, -1]]],
 [1, [[0, 1, 0], [0, -1, 0], [0, 0, 1]]],
 [1, [[0, 1, 0], [0, -1, 0], [0, 0, -1]]],
 [0, [[0, 1, 0], [0, -1, 0], [0, 0, 1], [0, 0, -1]]]]
sage: Arr_interior # The C_int regions
[[2, [[-1, 0, 1], [-1, 1, 0]]],
 [2, [[-1, 0, -1], [-1, 1, 0]]],
 [2, [[-1, -1, 0], [-1, 0, 1]]],
 [2, [[-1, -1, 0], [-1, 0, -1]]],
 [1, [[0, 0, 1], [0, 0, -1], [-1, 1, 0]]],
 [1, [[0, 0, 1], [0, 0, -1], [-1, -1, 0]]],
 [1, [[0, 1, 0], [0, -1, 0], [-1, 0, 1]]],
 [1, [[0, 1, 0], [0, -1, 0], [-1, 0, -1]]],
 [0, [[0, 0, 1], [0, 0, -1], [0, 1, 0], [0, -1, 0]]]]
\end{verbatim}
}
\end{example}

\begin{example}\rm  \label{ex:1F1}
Now let us provide a computation example of isomorphism classes for special functions.
The confluent hypergeometric function
$$
{}_1 F_1(a,c;x)= \sum_{k=0}^\infty \frac{(a)_k}{(1)_k (c)_k} x^k
$$
is annihilated by
\begin{equation}
L(a,c)=x \pd{x}^2 + (c-x) \pd{x} -a 
\end{equation}
It is obtained by restricting the GKZ hypergeometric system for
$A=\left(\begin{array}{ccc}
1&0&1\\
0&1&-1\\
\end{array}\right)$
to $x_1=x_2=1$ and by changing the variable $x_3 \mapsto -x$.
The correspondence with the parameter $\beta$ of the GKZ system is $\beta_1=-a, \beta_2=c-1$.
Put $M(a,c)=D/D L(a,c)$ where $D=D_1$.
Consider the start point $(a,c)$ and the target point $(a+1,c)$.
We have
\begin{eqnarray}
M(a,c)   &\stackrel{ x \pd{x}+a }{\longleftarrow}& M(a+1,c) \\
M(a+1,c) &\stackrel{ -x \pd{x}+x+a-c+1 }{\longleftarrow}& M(a,c) 
\end{eqnarray}
The composite of these left $D$-morphisms
\begin{equation*}
M(a,c) \ni \ell \mapsto \ell (-x\pd{x}+x+a-c+1) 
                \mapsto \ell (-x\pd{x}+x+a-c+1) (x\pd{x}+a) \in M(a,c)
\end{equation*}
is 
\begin{equation} \label{eq:conti-bf-1F1-a}
a(a-c+1),
\end{equation}
that is the $b$-function of this hypergeometric system for the direction $(1,0) \in \ZZ^2$.
Consider the direction $(0,-1)$.
We have
\begin{eqnarray}
M(a,c)   &\stackrel{ x \pd{x}+c-1 }{\longleftarrow}& M(a,c-1) \label{eq:1F1_conti_c1}\\
M(a,c-1) &\stackrel{ x \pd{x}-1 }{\longleftarrow}& M(a,c)  \label{eq:1F1_conti_c2}
\end{eqnarray}
and the $b$-function is  
\begin{equation} \label{eq:conti-bf-1F1-c}
a-c+1.
\end{equation}
We have four $2$-dimensional faces for the arrangement
$a (a-c+1)=0$ in $\QQ^2$.
Our algorithm classifies into isomorphism classes based on this arrangement.
\begin{enumerate}
\item \textbf{Exploration of highest dimensional chambers (\text{dim} = 2):} The algorithm extracts each 2-dimensional face (chamber) and evaluates the set of lattice points belonging to its strict interior, which excludes the boundary. For example, the region on the first quadrant side $C = \{a \ge 0, a-c+1 \ge 0\}$ becomes the target for exploration. 
\item \textbf{Shift by Hilbert basis and intersection check:}
To confirm whether all lattice points inside the chamber form an isomorphic class, shifts along the Hilbert basis of the interior are verified. For instance, for the above $C$, we calculate the contiguity relation for the shift $h=(0,-1)$, i.e., $(a,c) \to (a,c-1)$. From the $b$-ideal corresponding to this shift, a candidate for a new wall that breaks isomorphism, $W = \{-a+c-1=0\}$, is derived. 
    
    The algorithm then determines whether this $W$ crosses the interior of the current chamber ($\textsc{Intersection}(W, C_{int}) \neq \emptyset$). In this case, since it is confirmed that they do not intersect, the chamber is not further subdivided, and all points inside are determined to be mutually isomorphic (the vast red, yellow and green yellow regions in Figure \ref{fig:c1f1}).
\end{enumerate}

Once the exploration of the 2-dimensional regions is complete, the algorithm recursively transitions to investigating the boundaries of the chambers (codim 1 faces).
We consider an arrangement
on the $1$-dimensional space $L=\{a=0\}$. 
The contiguity relation on $V$ is also given
by (\ref{eq:1F1_conti_c1}) and
(\ref{eq:1F1_conti_c2}).
In this case as well, no wall points appear other than those originating from the initial walls.
Then, the arrangement has two $1$-dimensional faces and one $0$-dimensional face $(0,1)$.

We also consider arrangement
on $L=\{a-c+1=0\}$.
The contiguity relation is given as 
\begin{eqnarray}
M(c-1,c)   &\stackrel{ \pd{x} }{\longleftarrow}& M(c,c+1) \\
M(c,c+1) &\stackrel{ x \pd{x}-x+c }{\longleftarrow}& M(c-1,c) 
\end{eqnarray}
and the $b$-function is $c-1$.
Then, the arrangement has two $1$-dimensional
faces and two $0$-dimensional faces
$(0,1)$.
In other words, in this case too, no wall points appear other than those derived from the initial walls.

Finally, intersections of codim 2 (e.g., $a=0, c=1$) are evaluated, completing the classification across all dimensions.

Figure \ref{fig:directed-graph-1F1} illustrates isomorphic classes
by our method.
Big circles of the figure are reduced set of representatives.
Figure \ref{fig:c1f1} visualizes the output results generated by our implementation, where the lattice points inside regions (as well as lines and points) painted in the same color correspond to mutually isomorphic $D$-modules.

\begin{figure}[tb]
\begin{center}
\begin{tikzpicture}[scale=0.6]
\draw (0,5) node[right]{$a=0$};
\draw[thick,black] (0,-5)--(0,6);

\draw (4,5) node[right]{$c=a+1$};
\draw[thick,black] (-5,-4)--(5,6);

\draw(2,1) node[above]{$c=1$};
\draw[dotted] (-2,1)--(2,1);

\filldraw [pattern=horizontal lines gray] (-1,-1) circle (0.1);
\filldraw [pattern=horizontal lines gray] (1,-1) circle (0.1);

\filldraw [fill=white] (-1,0) circle (0.1);
\filldraw [fill=black] (0,0) circle (0.1);

\filldraw [pattern=horizontal lines gray] (-1,1) circle (0.1);
\filldraw [fill=white] (0,1) circle (0.1);

\filldraw [fill=black] (0,2) circle (0.1);
\filldraw [fill=white] (1,2) circle (0.1);

\filldraw [pattern=horizontal lines gray] (1,3) circle (0.1);

\filldraw [pattern=horizontal lines light gray] (1,-1) circle (0.2);
\filldraw [pattern=horizontal lines light gray] (0,0) circle (0.2);
\filldraw [pattern=horizontal lines light gray] (-1,1) circle (0.2);
\filldraw [pattern=horizontal lines light gray] (1,2) circle (0.2);

\draw[->] (-1,-1)--(0-0.1,-1);

\draw[->] (0,0)--(0,-1+0.1);

\draw[->] (-1,1)--(0-0.1,1);
\draw[->] (-1,1)--(-1,0+0.1);

\draw[->] (-1,2)--(-1,1+0.1);
\draw[->] (-1,2)--(0-0.1,2);

\draw[->] (1,3)--(1,2-0.1);

\draw[->] (1,2)--(2-0.1,3-0.1);
\draw[->] (-1,0)--(-2,-1);
\end{tikzpicture}
\end{center}
\caption{Representatives for ${}_1F_1$}
\label{fig:directed-graph-1F1}
\end{figure}

\begin{figure}[htbp]
\begin{minipage}[b]{0.48\textwidth}
    \centering
    \includegraphics[width=0.6\linewidth]{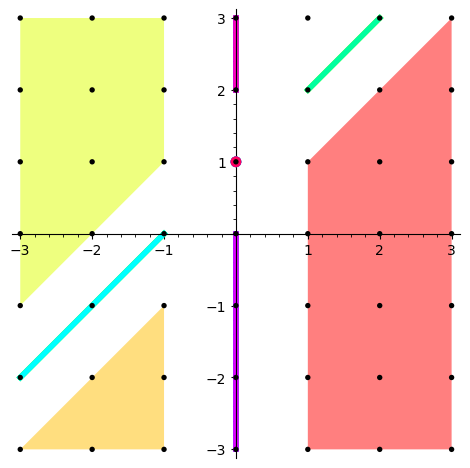}
    \captionof{figure}{Isomorphism classes of ${}_1F_1(a,c;-x)$ in the parameter space $(a,c)$. Distinct colors and line segments represent different isomorphism classes separated by the walls.}
    \label{fig:c1f1}
\end{minipage}
\hfill
\begin{minipage}[b]{0.48\textwidth}
    \centering
    \includegraphics[width=0.6\linewidth]{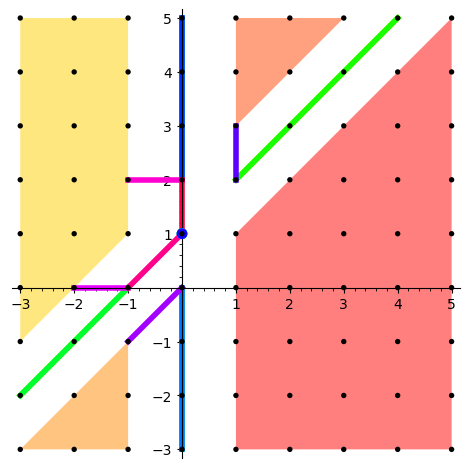}
    \captionof{figure}{Applying merging algorithm to isomorphism classes of ${}_1F_1(a,c;-x)$ in the parameter space $(a,c)$ obtained by \textsc{Representatives}}
    \label{fig:merged_1F1_isom}
\end{minipage}
\end{figure}
\end{example}

\begin{proof}[Proof of Theorem \ref{th:correctness}]
Horn type hypergeometric systems are defined as those obtained by restricting GKZ hypergeometric systems to $x_1=\cdots=x_d=1$.
Therefore, we first note that if the GKZ hypergeometric systems are isomorphic, the corresponding Horn type hypergeometric systems are also isomorphic.

If the algorithm terminates, every face of the obtained arrangement is either bounded, or its interior points are isomorphic. Therefore, (1) holds.

Theorem 5.2 in \cite{saito-isom} shows that the condition for isomorphism of GKZ systems can be entirely described by inequality conditions concerning primitive supporting functions.
Therefore, walls that cannot be written via primitive supporting functions will never appear, and Algorithm \ref{alg:rep} will not be recursively called by Algorithm \ref{alg:rep_sub} indefinitely. Thus, (2) holds.
\end{proof}

\subsection{Addition of Walls}  \label{sec:addition_of_walls}
The function \textsc{IsAffineSpaceIncludedInArrangement}(W,L,E) considers $V_L$, which is the zero set of the arrangement in the ambient space defined by $E$ restricted to the affine subspace $L$, and determines whether $V_L$ includes $W$.
If it does not, a new wall is constructed using the function $E'$=\textsc{GetNewHyperplanes(W,L,E)}.
In other words, this constructs a wall such that $V'_L$, the restriction of the zero set of the arrangement defined by $E' \cup E$ to the affine subspace $L$, includes $W$.
Note that this method of construction is not unique.

\medbreak

Let us add an explanation using an example to show the behavior when the codimension of a wall that might break isomorphism is 2 or more.
\begin{example}\rm \label{ex:f0134-a}
This is an example of a 2-variable hypergeometric system restricted to $x_1=x_2=1$ from the GKZ system $H_A(\beta)$ defined by 
$A=\left(\begin{array}{cccc} 1&0&-2&-3\\ 0&1&3&4\\ \end{array}\right)$. 
Here $I_A$ is not Cohen-Macaulay.
The so-called hole is located at $(-1,2)$.
The parameters are $(b_1,b_2)$.
Let $C$ be the 2-dimensional face $-2+b_2 \geq 0, 1+b_1 \geq 0$.
Then $C_{int}$ is $-3+b_2 \geq 0, b_1 \geq 0$.
In our implementation, the internal representation of $C$ is {\tt [2, [[-2,0,1],[1,1,0]]]}.
That means we represent $c_0+c_1 b_1 + c_2 b_2 \geq 0$ as {\tt $[c_0,c_1,c_2]$}.
Now, if we compute \textsc{Isom}(param=[b1,b2],start=[b1,b2],target=[b1,b2+1]) in the direction of $(0,1)$ to find the walls, we get
{\tt 
W=[[-3, -4, -3]], W=[[-2, -4, -3]], W=[[-1, -4, -3]], W=[[1, 0, 1]], W=[[-2, 0, 1], [1, 1, 0]]}. 
Here, for example, {\tt [-3,-4,-3]} signifies the plane $-3-4b_1-3b_2=0$.
What requires attention is the final entry {\tt W=[[-2, 0, 1], [1, 1, 0]]}.
This represents the point $(-1,2)$ (the hole).
In this case, since the intersection of the interior $C_{int}$ of $C$ and these $W$ is empty, updating the set of hyperplanes $E$ to $E'$ does not occur. However, if $(-1,2)$ were contained in $C_{int}$, then $E'$ would be supplemented with, for example, $-2 + b_2=0$ and $4b_1+3b_2-2=0$.
Here $b_2=0, 4b_1+3b_2=0$ are primitive supporting hyperplanes.
\end{example}

A point of caution when adding walls is that the representation of $(-1,2)$ as an intersection of hyperplanes is not unique.
If a wall is represented as an intersection of hyperplanes, it should be implemented such that it is represented uniquely.
Similar consideration is needed for representing the affine subspace $L$.
Furthermore, the same care is required for the intersection between $L$ and $W$.

One method to address this is to use the representation method for walls from Theorem 4.5 and Proposition 4.8 of \cite{saito-isom}.
In this representation, a wall is contained within the union of translations of zero sets of primitive supporting hyperplanes.
This method is detailed further in subsection \ref{sec:b-chi}.

By using this heuristic, the proliferation of the new hyperplane arrangement $E'$ can be kept to a minimum.

According to Theorem 2.1 in \cite{saito-isom}, the number of isomorphism classes is finite.
Therefore, it is expected that there is an algorithm that terminates in finitely many steps, but we have not yet obtained a proof that our algorithm with these heuristics will terminate in a finite number of steps.

\subsection{Algorithm for Computing $B_{\chi,-\chi}$} \label{sec:b-chi}

M. Saito \cite[Th 4.9]{saito-isom} introduced the monomial ideal in $\QQ[\pd{1}, \ldots, \pd{n}]$ defined by
\begin{equation}
M_\chi = \langle \pd{}^u \,|\, Au \in \chi + \NN_0 A \rangle
\end{equation}
for $\chi \in \ZZ A$, and using this monomial ideal, gave a sufficient condition for two GKZ systems $M_A(\beta)$ and $M_A(\beta+\chi)$ to be isomorphic.

In this section, we provide an algorithm to compute $M_\chi$.
This computation reduces to finding the entire set of monomials $J_M$ contained in the ideal $J$ generated by the toric ideal $I_A$ and the monomial $\pd{}^m$, proceeding as follows:
\begin{enumerate}
\item Find a pair $U,V \in \NN_0^n$ such that $\chi=AU-AV$.
\item By definition, there exists a $W \in \NN_0^n$ such that $Au = AU-AV+AW$.
\item Because $A(u+V)-A(U+W)=0$, it follows that $\pd{}^{u+V}-\pd{}^{W+U} \in I_A$.
  Therefore, when $m=U$,
$$ \pd{}^{u+V} \in \langle I_A, \pd{}^U \rangle =J$$  
\item The ideal quotient $J_M:\pd{}^V$ of the monomial ideal $J_M$ contained in $J$ provides the desired $M_\chi$.
\end{enumerate}

\begin{algorithm}[H]
\caption{\textsc{monomialPart}}
\label{alg:monomialPart}
\begin{algorithmic}[1]
\Require $A$ to define a toric ideal, monomial $\pd{}^m$.
\Ensure Monomial part $J_M$ of $J=\langle I_A, \pd{}^m\rangle$.
\State Fix an order $\succ$ and let $G$ a Gr\"obner basis of $J$ by the order.
\State $G=G_M \cup G_B$ where $G_M$ is a set of monomials in $G$ and $G_B$ is a set of binomials in $G$.
\State $M'=G_M$, $M=\emptyset$.
\While{ $M' \not= M$}
   \State $M=M'$.
   \State $M'=M \cup \cup_{\pd{}^u-\pd{}^v \in G_B, \pd{}^u \succ \pd{}^v} (M:\pd{}^v)\pd{}^u$
\EndWhile
\State \Return $M$
\end{algorithmic}
\end{algorithm}

\begin{example}\rm
For $A=\left(\begin{array}{ccc} 2&1&1\end{array}\right)$,
we have $I_A=\langle \pd{2}-\pd{3}, \pd{3}^2-\pd{1} \rangle$.
\begin{enumerate}
\item From the Gr\"obner basis of $\langle I_A, \pd{1} \rangle$ with respect to the graded reverse lexicographic order, we get
$$ G_M=\{\pd{1}, \pd{3}^2\}, \ G_B=\{\underline{\pd{2}}-\pd{3} \} $$
where the underline indicates the leading term.
Thus $M=\{ \pd{1}, \pd{3}^2 \}$.
\item Since $ \pd{3}^2:\pd{3}=\pd{3} $ and $\underline{\pd{2}}-\pd{3} \in G_B$, we obtain $M'=M \cup \{ \pd{3} \pd{2} \}$.
\item The new $M$ is $\{ \pd{1}, \pd{3}^2, \pd{2}\pd{3} \}$.
\item Since $\pd{2}\pd{3}:\pd{3}=\pd{2}$, we obtain
$M'=M \cup \{ \pd{2} \pd{2} \}$.
\item The ideal quotient does not increase the set any further, so the returned value is
$\{ \pd{1}, \pd{3}^2, \pd{2}\pd{3}, \pd{2}^2 \}$.
\end{enumerate}
\end{example}

\begin{theorem} \label{th:correctness_monomial_part}
The output of Algorithm \textsc{monomialPart} is the set of generators for the ideal $J_M$ formed by the monomials contained in $\langle \pd{}^m, I_A \rangle$.
\end{theorem}

\begin{proof}
The proof relies on "reverse generation of reduction."
Assume $\pd{}^e \in J_M$.
$\pd{}^e$ can be reduced to $0$ using $G_M$ and $G_B$.
When reduced by a binomial belonging to $G_B$, the monomial is reduced to another monomial.
When reduced by a monomial belonging to $G_M$, the monomial is reduced to $0$.
From the process of this reduction, we have
\begin{equation} \label{eq:toric_reduction}
  0
= \pd{}^e - \pd{}^{w^{(1)}}(g_1^{(1)}-g_2^{(1)})
          - \pd{}^{w^{(2)}}(g_1^{(2)}-g_2^{(2)})
          - \cdots
          - \pd{}^{w^{(s)}}(g_1^{(s)}-g_2^{(s)}) - \pd{}^{w^{(s+1)}} h^{(s+1)}
\end{equation}
Here, $g_1^{(i)}-g_2^{(i)}$ is an element of $G_B$, and $g_1^{(i)} \succ g_2^{(i)}$.
Also, $h^{(s+1)}$ is an element of $G_M$.
Let the sets $M$ appearing in each step of Algorithm \ref{alg:monomialPart} be $M_0, M_1, M_2, \ldots$ (with initial value $M_0 = G_M$).
From the last relation of equation (\ref{eq:toric_reduction}), we have $\pd{}^{w^{(s)}} g_2^{(s)} = \pd{}^{w^{(s+1)}} h^{(s+1)}$, which belongs to the ideal generated by $M_0$.
Therefore, by the algorithm's definition, the pre-reduction term $\pd{}^{w^{(s)}} g_1^{(s)}$ belongs to the ideal generated by $M_1$.
Here, $\pd{}^{w^{(s)}} g_1^{(s)}$ is equal to the result of the previous reduction $\pd{}^{w^{(s-1)}} g_2^{(s-1)}$, so this also belongs to the ideal generated by $M_1$.
We continue to decrement $s$ in this manner. That is, we can show inductively that $\pd{}^{w^{(s-k)}} g_2^{(s-k)}$ belongs to the ideal generated by $M_k$, and the pre-reduction term $\pd{}^{w^{(s-k)}} g_1^{(s-k)}$ belongs to the ideal generated by $M_{k+1}$.
Since $J_M$ is finitely generated, the possible values for $s$ are also finite.
Consequently, the algorithm outputs the generators of $J_M$.
\end{proof}

Let $F_\sigma(s)$, $s\in \QQ^d$ be the primitive supporting hyperplane for the facet $\sigma$ of the cone $\QQ_{\geq 0} A$.
Let $S(M_\chi)$ be the set of standard pairs of the monomial ideal $M_\chi$ (see, e.g., \cite{SST}).
\begin{theorem}[{\cite[Th 4.5, Th 4.9]{saito-isom}}]
When $\chi \in \ZZ A$ and $\beta \not \in V(B_{-\chi,\chi})$, the two GKZ systems $M_A(\beta)$ and $M_A(\beta+\chi)$ are isomorphic as left $D$-modules.
Here,
\begin{equation} \label{eq:V_B}
V(B_{-\chi,\chi})=(V(B_\chi)-\chi) \cup V(B_{-\chi})
\end{equation}
\begin{equation} \label{eq:V_B_F}
V(B_\chi) = \cup_{(u,\tau) \in S(M_\chi)} 
   V\left(\langle F_\sigma-F_\sigma(Au) \,|\, 
          \mbox{$\sigma$ facet} \supset \tau \rangle\right)
\end{equation}
\end{theorem}

In \textsc{GetNewHyperplanes} of Algorithm \ref{alg:rep}, our heuristic is to find the affine variety (\ref{eq:V_B}) containing the wall $W$, and to add $F_\sigma(s)-F_\sigma(Au)=0$ from (\ref{eq:V_B_F}) to the hyperplane set $E'$.
We conjecture that this method guarantees the termination of our algorithm, but we have not yet proven this.

\subsection{Merging Isomorphism Classes}  \label{sec:merging}

The isomorphism classes determined by the \textsc{Representatives} algorithm may be registered as distinct entities even if they belong to the same isomorphism class.
The \textsc{MergeIsomClasses} algorithm discovers these identical classes via exhaustive search and finds the isomorphisms connecting them.

\begin{algorithm}[H]
\footnotesize
\caption{\textsc{MergeIsomClasses}}
\label{alg:merge}
\begin{algorithmic}[1]
\Require $IsomClass, Contiguity, IsomClassClosure$.
\Ensure $MergedIsomClass, MergedContiguity$
\State $MergedIsomClass \gets \textsc{SimplifyFaces}(IsomClass)$
\State $MergedContiguity \gets Contiguity$
\For{\textbf{each} pair of indices $(i, j)$ in $IsomClass$}
    \If{$\textsc{IsContainedWithCodim1}(IsomClassClosure[i], IsomClass[j])$}
        \State $p \gets \text{a base point of } IsomClass[j]$
        \State $q \gets \textsc{FindNearbyPoint}(p, IsomClass[i])$
        \If{$q \neq \text{null}$}
            \State $I_{b} \gets \textsc{Isom}([],p, q)$ \Comment{Evaluate $b$-ideal at points $p, q$}
            \If{$I_{b} \neq \langle 0 \rangle$} \Comment{$p$ and $q$ belong to the same isomorphism class}
                \State $seg \gets \textsc{SegmentToFace}(p, q)$
                \If{$seg \notin MergedIsomClass$}
                    \State Append $seg$ to $MergedIsomClass$
                    \State $C_{pq} \gets \textsc{Contiguity}(p, q)$
                    \State $C_{qp} \gets \textsc{Contiguity}(q, p)$
                    \State Append $(C_{pq}, C_{qp}, q-p)$ to $MergedContiguity$
                \EndIf
            \Else
               \State Call Tsai-Walther's algorithm \cite{tw-hom} 
               \State \quad to check if two $D$-modules are isomorphic.
               \State If they are isomorphic, append $seg$ to $MergedIsomClass$ and add contiguity.
            \EndIf
        \EndIf
    \EndIf
\EndFor
\State \Return $(MergedIsomClass, MergedContiguity)$
\end{algorithmic}
\normalsize
\end{algorithm}

An efficient implementation of Tsai-Walther's algorithm remains a task for the future, and is currently not included in our implementation.

\begin{example} \rm
Applying this algorithm to ${}_1F_1$ yields the result shown in Figure \ref{fig:merged_1F1_isom}. Consequently, there are 4 isomorphism classes.

\end{example}

\begin{figure}[htbp]
    \centering
    \includegraphics[width=0.7\linewidth]{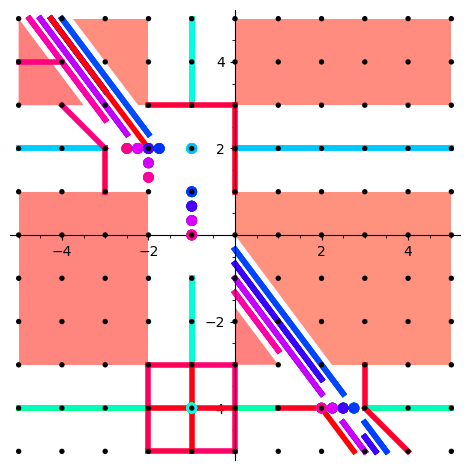}
    \caption{Merged isomorphism classes for the restricted 2-variable hypergeometric system derived from the GKZ system with a hole at $(-1,2)$. Bold segments indicate contiguous isomorphic points across codim 1 boundaries, successfully connected by Algorithm~\ref{alg:merge}.}
    \label{fig:f0134-merged}
\end{figure}

\begin{example}\rm 
(Continuation of Example \ref{ex:f0134-a}.)

Figure \ref{fig:f0134-merged} displays the isomorphism classes of a 2-variable hypergeometric system restricted from a GKZ hypergeometric system that is given by the configuration matrix 
$$ A = \begin{pmatrix} 1 & 0 & -2 & -3 \\ 0 & 1 & 3 & 4 \end{pmatrix} $$
and possesses a hole at $(-1,2)$ in its parameter space. 

The initial isomorphism classes (chambers) extracted by Algorithm~\ref{alg:rep} are geometrically divided by walls of codimension 1. However, algebraically, isomorphic modules may exist across these walls. Algorithm~\ref{alg:merge} (\textsc{MergeIsomClasses}) explores pairs of lattice points $(p, q)$ crossing the boundary (codim 1 face) between adjacent chambers and determines their isomorphism by evaluating the $b$-ideal. 

If they are confirmed to be isomorphic, the algorithm generates a line segment (a 1-dimensional face) connecting these points and merges the spaces. The bold segments (edges) in Figure~\ref{fig:f0134-merged} represent the bridging contiguity relations discovered by this merging process, visually proving the non-trivial connected structures of the isomorphism classes surrounding the hole. 
\end{example}

\section{Comprehensive Restriction Algorithm} \label{sec:restriction}

Let $D_n$ be the Weyl algebra of $n$ variables.
Let $I(\kappa)$ be a holonomic left ideal of $D_n$ 
with parameters $\kappa \in \CC^d$.
We want to compute the restriction module
\begin{equation}  \label{eq:restriction_module}
 \frac{D_n}{I(\kappa)+x_1 D_n + \cdots + x_m D_n}.
\end{equation}

\begin{algorithm}[H]\rm 
\footnotesize
\caption{Comprehensive restriction algorithm that gives 
a partial answer}   \label{alg:alg-rest-top} \rm
\begin{itemize}
\item {\bf Require}: $I(\kappa)$, $x_1=\cdots=x_m=0$. 
\item {\bf Ensure}: Strata $S_1$, $S_2$ and $S_3$ of the $\kappa$ space $\CC^d$.
The restriction module (\ref{eq:restriction_module}) on each stratum of them. 
\end{itemize}
\begin{enumerate}
\item {\tt Ans=[\,]}.
\item Put $w=(\overbrace{1,\ldots,1}^{m}, 0, \ldots, 0)$ and compute
a comprehensive Gr\"obner system $G$ by \, $\prec_{(-w,w)}$ order.
\item Compute comprehensive $b$-functions for restriction, which are monic generators
of \, ${\rm in}_{(-w,w)}(I(\kappa)) \cap \CC[\theta_1+\cdots+\theta_m]$.
Let $S_1$ be strata of the comprehensive $b$-functions that refines strata 
of the comprehensive Gr\"obner system.
\item For each stratum $U$ of $S_1$, refine $U$ into subsets such that 
(a) the maximal non-negative integral root
of the $b$-function is $0$ or (b) no non-negative integral root of $b$ or (c) other cases on each subset.
Let $S_2'$ be the collection of subsets such that (a) or (b) holds.
Let $S_3'$ be the collection of subsets such that (c) holds.
\item For each stratum $V$ of $S_2'$,
compute a comprehensive Gr\"obner system of $G''=G'|_{x_1=\cdots=x_m=0}$
where $G'$ is the collection of the elements of the $(-w,w)$ Gr\"obner basis $G$
on $V$
such that ${\rm ord}_{(-w,w)}(g) \leq 0$,  $g \in G$.
It refines $V$ and let $S_2$ be the collection of these refinement.
\begin{enumerate}
\item The restriction module on each stratum $W$ of $S_2$ of type (a) is
$\frac{D'}{\left\langle \mbox{\rm (a Gr\"obner basis of $G''$ on $W$)} \right\rangle}$
where $D'=\CC\langle x_{m+1}, \ldots, x_n, \pd{m+1}, \ldots, \pd{n} \rangle$.
Append them to {\tt Ans}.
\item The restriction module (output) is $0$ on the set of strata of $S_2'$ of type (b). Append them to {\tt Ans}.
\end{enumerate}
\item For $V$ in $S_3'$ (type (c)), 
if $I(\kappa)$, $\kappa \in V$ is a GKZ system or a hypergeometric system of Horn type discussed in previous sections, 
call 
\, {\tt Rfr}($I(\kappa(p)), p,G, V)$ (\underline{r}epresentatives \underline{f}or \underline{r}estriction, Algorithm \ref{alg:hg-rest})
where we reparametrize $\kappa$ by $p$ as $\kappa_i(p)=p_i$.
The return value is strata $S_3$ and restrictions on each stratum of $S_3$.
Append them to {\tt Ans}.
If $I(\kappa)$ does not belong to these hypergeometric systems, this algorithm
does not give an answer for type (c) strata $S_3'$.
\item Return {\tt Ans}.
\end{enumerate}
\normalsize
\end{algorithm}

\begin{example}\rm
Let $\kappa=(a,b,c) \in \CC^3$
and consider the left ideal $I(\kappa)$ of $D_1$
generated by the Gauss hypergeometric operator
$ L=x(1-x)\p^2 + (c - (a + b + 1)x)\p - ab $.
Here $x_1$ is denoted by $x$ and $\pd{1}$ by $\p$.
The set $\{L\}$ is $(-1,1)$ Gr\"obner basis for any $\kappa$.
Then $S_1 = \{ \CC^3 \}$.
We have ${\rm in}_{(-1,1)} (L)= x \p^2 + c \p$,
then the $b$-function for restriction is $\theta_1 (\theta_1+c-1)$.
Then $S_2'$ is $\{ V \}$ where
$V=\{ (a,b,c) \,|\, c \not\in \ZZ_{\leq 0} \}$.
Since ${\rm ord}_{(-1,1)}(L)=1$, we have $G'=\emptyset$.
Then the restriction module on $V$ is $D'=\CC$.
The strata $S_3'$ (case (c)) is $\{W\}$ where 
$W=\{ (a,b,c) \,|\, c \in \ZZ_{\leq 0} \}$.
This case will be discussed in Section \ref{sec:gauss-hg}.
\end{example}

Let ${\bar q}_1$ be the maximal non-negative integral root
of the $b$-function for the restriction.
Let $G$ be the $G$ that appears in Algorithm \ref{alg:alg-rest-top}.
We denote by 
$\mbox{\tt Rest}(G,{\bar q}_1)$
the output by CGS of the final step of 
computing Gr\"obner basis in a free module
of the restriction algorithm, see, e.g.,
\cite[Steps 6 and 7 of Alg. 5.2.8]{SST}.

\begin{example} \rm
Let $H_2(a,b,b',c,c')$ be the left ideal
generated by (\ref{eq:f2a}) and (\ref{eq:f2b})
that annihilates the Appell function $F_2$.
Consider the left $D_2$-module $M(a,b,b',c,c')=D_2/H_2(a,b,b',c,c')$.
Suppose that $c=0$ and $c' \notin \ZZ_{\leq 0}$.
The maximum non-negative root of the $b$-function is $s_0=1-c=1$
(see Example \ref{ex:f1-f2-comprehensive-b-fn}).
The restriction module is
$$
\frac{\CC\pd{x}+\CC\pd{y}+\CC}
 {
   \CC (-a b) + \CC c' \pd{y} + \CC (-a b')
 }.
$$
The dimension is equal to $3$ minus the rank of the matrix
$$
\begin{pmatrix}
0 & 0 & -\aalpha \bbeta \\
0 & \ggamma' & -\aalpha \bbeta' 
\end{pmatrix}.
$$
The stratification with respect to the rank can be obtained
by a comprehensive Gr\"obner system for linear polynomials.
$\mbox{\tt Rest}(G,1)$ returns this comprehensive Gr\"obner basis.
\end{example}

In order to compute the restriction modules on the strata $S_3'$
for a GKZ system or for a hypergeometric system of Horn type, 
we apply the algorithm \ref{alg:hg-rest} utilizing algorithms to find contiguity relations (isomorphisms).
This algorithm is a variation of \textsc{Representatives}
(Algorithm \ref{alg:rep}).

\begin{algorithm}[H]  \rm
\footnotesize
\caption{Comprehensive restriction: procedure $\mbox{\tt Rfr}(H_A(\beta(p),p,G',E)$} \label{alg:hg-rest} 
{\bf Require}: $H_A(\beta(p))$ (hypergeometric system), $p$ (a set of $m$ parameters), 
$G'$ ($(-w,w)$-Gr\"obner basis), $E$ (conditions). \\
{\bf Ensure}: a list of [conditions(stratum), restriction, contiguity relations].
\begin{enumerate}
\item {\tt Ans=[\,]}.
\item Let $r_1(p), \ldots, r_k(p)$ be the roots
of the $b$-function for the restriction.
\item For all $r_i$, assume $r_i \in \ZZ_{\geq 0}$ or not
and relationships (larger or smaller or equal)
among $r_i$'s supposed to be non-negative integers.
Let $\tilde K$ be the set of all distinct assumptions
on $r_i$'s.
\item For $K \in {\tilde K}$ do
\begin{enumerate}
\item $q(p)$ be the maximal non-negative integral root under the assumption $K$.
If $q(p)$ is a constant, append
$[E \cap K, \mbox{\tt Rest}(G,q(p)), \emptyset]$ to {\tt Ans}
and continue the for-loop.
\item Changing the indexing, 
we suppose that $q(p)$ depends on $p_1$.
Introduce new variables $q_1, \ldots, q_m$
such that
$q_1 = q = \alpha p_1 + \cdots$ ($\alpha \not= 0$),
$q_i=p_i$ ($i \geq 2$)
and express $\beta$ by $q$.
\item Choose $\delta \in \ZZ_{>0}$ so that
\begin{equation}
 \beta_i(q+\delta e_1)-\beta_i(q) \in \ZZ \quad
\mbox{for all $i$}
\end{equation}
where $e_1=(1,0, \ldots, 0)$.
\item Put $\Lambda = \{0, 1, 2, \ldots, \delta-1\}$.
\item For $k$ in $\Lambda$ do
\begin{enumerate}
\item Derive the contiguity relation for it $(L_u, L_d, b(q))$
\begin{equation}
  D/H_A(\beta(q+k e_1))
\stackrel{L_d}{\rightarrow}
  D/H_A(\beta(q+(k+\delta)e_1))
\stackrel{L_u}{\rightarrow}
  D/H_A(\beta(q+k e_1))
\end{equation}
where $ L_d L_u \equiv b(q)$.
\item Consider 
\begin{equation}  \label{eq:segments}
 \left( \RR^m \setminus V(b) \right) \cap
 \left( \RR \times ({\bar q}_2, \ldots, {\bar q}_m) \right)
\end{equation}
where ${\bar q}_2, \ldots, {\bar q}_m$ are generic numbers.
Let $Q$ be the set of the $q_1 \in \NN_0$ that is the minimum 
in each first coordinate of connected components
of (\ref{eq:segments}).
\item Factorize $b(q)$ into degree $1$ polynomials as $\prod_{j=1}^J b_j(q)$.
\item For all ${\bar q}_1 \in Q$ do
\begin{enumerate}
\item Append $[E \cap K \cap \{ b(q)\not=0 \}, \mbox{\tt Rest}(G,{\bar q}_1),\mbox{the contiguity relation}]$.
\end{enumerate} 
\item For all factors $b_j$ in $b$ do
\begin{enumerate}
\item Eliminate one variable in $q_1, \ldots, q_m$ by $b_i(q)=0$.
Changing indices, we suppose that $q_m$ is eliminated.
\item Express $\beta$ in terms of $q_1, \ldots, q_{m-1}$.
\item Append \, $\mbox{\tt Rfr}(H_A(\beta(q_1, \ldots, q_{m-1}),
(q_1, \ldots, q_{m-1}), G, E \cap K \cap \{b_i=0\})$
to {\tt Ans}.
\end{enumerate}
\end{enumerate}
\end{enumerate}
\item Return {\tt Ans}.
\end{enumerate}
\normalsize
\end{algorithm}

\begin{remark}\rm
Although, as long as we have tried, we can always find a contiguity
relation of $\delta e_1$ shift,
we might fail at this step.
If we fail to find a contiguity relation of $\delta e_1$ shift,
we need to increase $\delta$.
Since the number of isomorphic classes of $H_A(\beta+\iota)$, $\iota \in \ZZ^d$
are finite by \cite{saito-isom}, we can find a contiguity relation
at a suitable $\delta e_1$.
\end{remark}

\begin{remark}\rm
Algorithms \ref{alg:alg-rest-top} and \ref{alg:hg-rest} 
will be generalized to obtain
a restriction complex (a restriction of $D_n/I(\kappa)$ in a derived category)
by applying the algorithm of \cite{OT-rest}.
Note that we need to replace ``maximal non-negative integral root''
of the algorithms \ref{alg:alg-rest-top} and \ref{alg:hg-rest}
by ``maximal integral root''.
A comprehensive version of $(-w,w)$-adapted resolution is an open question
to give an algorithm to obtain a restriction complex.
\end{remark}

\begin{example}\rm \label{ex:1F1-restriction}
This is a continuation of Example \ref{ex:1F1} (${}_1 F_1$ case).
The $b$-function for the restction to $x=0$
is $s(s+c-1)$.
The roots are $s=0$ and $s=1-c$.
Type (a) case is $1-c \not\in \ZZ_{\geq 0}$.
Since ${\rm ord}_{(-1,1)}(L)=1$, $G''$ is empty.
Then the restriction is $\CC$.
Since $s$ is a factor of the $b$-function, type (b) case does not occur.

Consider the type (c) case.
In other words, assume $c \in \ZZ_{\leq 1}$.
Let this assumption be $E$.
$G$ is $\{ L \}$.
We call the procedure 
$\mbox{\tt Rfr}(H_A((a,c)), (a,c), G, E)$.
Put $q_1=1-c$ and $q_2=a$.
Then, $\Lambda=\{ 0 \}$.
The condition on $c$ becomes $q_1=1-c \in \ZZ_{\geq 0}$
and the roots of the $b$-function for restriction is $0$ and $q_1$.
Firstly, we compute 
a contiguity relation for the shift from $q_1=1-c$ to $q_1+1=1-(c-1)$.
The $b$-function for contiguity is $q_1+q_2=1-c+a$
by (\ref{eq:conti-bf-1F1-c}).
SInce we assume ${\bar q}_2 \not \in \ZZ$,
$\left( \RR^2 \setminus V(q_1+q_2) \right) \cap \left( \RR \times {\bar q}_2 \right)$ 
is $ (-\infty,\infty) \times {\bar q}_2$.
Then the set $Q$ of the minimal non-negative integers in the connected component is $\{ 0 \}$.
The restriction on this stratum is isomorphic to that of $D/H_A((a,c)=(a,0))$.
Since the maximal integral root is $1$ in the $(a,c)$ parameter space
and ${\rm ord}_{(-1,1)}(L)=1$,
we have $G''=\{ L|_{x=0}= -a \}$.
Thus, the restriction is $\CC$ when $a \not= 0$
and is $\CC^2$ when $a=0$.
Secondly, we consider the case $q_1+q_2 = 1-c+a=0$.
The parameter space is one dimensional and parametrized as
$(a,c)=(0,1)-(1,1)s'$.
We call the procedure
$\mbox{\tt Rfr}(H_A((-s',1-s')), s', G, \{ c \in \ZZ_{\leq 0}, c=a+1=-s'+1 \}$.
The $b$-function for the contiguity
of the shift $s' \mapsto s'+1 \mapsto s'$ 
is $s'(s'-1)=c(c-1)$.
Then, the cases of $s'=\{ -1,0,1 \}$ are representatives
of isomorphic classes.
In other words, $(a,c)=(1,2), (0,1), (-1,0)$ are the representatives. 
The restrictions are all $\CC$. 
\end{example}

This procedure will be a little complicated.
Then, more examples will help.
The comprehensive restriction algorithm will be illustrated for
the Gauss hypergeometric system and the system of Appell function 
$F_1$ 
in the following sections \ref{sec:gauss-hg} and \ref{sec:f1}.

\section{Restriction of the Gauss Hypergeometric System to the Origin} \label{sec:gauss-hg}

The Gauss hypergeometric function ${}_2F_1(a,b,c;x)$ is annihilated by
the operator
\begin{equation} \label{eq:gauss_hg_op}
L(\aalpha, \bbeta, \ggamma) = x(1-x)\pd{x}^2 + (\ggamma - (\aalpha + \bbeta + 1)x)\pd{x} - \aalpha \bbeta.
\end{equation}
We consider the left ideal generated by $L$
$$
\Hg(\aalpha, \bbeta, \ggamma) = D L(\aalpha, \bbeta, \ggamma)
$$
where $D=D_1$.
We will compute the restriction module
$$M(\aalpha, \bbeta, \ggamma)/xM(\aalpha, \bbeta, \ggamma) \cong D/(\Hg(\aalpha, \bbeta, \ggamma) + x D)$$ 
of the left $D$-module $M(\aalpha, \bbeta, \ggamma) = D/\Hg(\aalpha, \bbeta, \ggamma)$
to $x = 0$.

The generic $b$-function (for restriction)
is 
$$
b(s) = s(s+\ggamma-1)
$$
with respect to the weight vector $w=(1)$.
The stratum for this $b$-function is 
$\CC^3=\{ (\aalpha, \bbeta, \ggamma) \in \CC^3 \}$.
The maximal non-negative integral root $s_0$ of
$b(s)$ is
$$
s_0 = 
\begin{cases}
0 & (\ggamma \notin \mathbb{Z}_{\leq 0}) \\
1-\ggamma & (\ggamma \in \mathbb{Z}_{\leq 0}).
\end{cases}
$$

The Gr\"obner basis $G$ of $\Hg$ by Algorithm \ref{alg:alg-rest-top}
is $\{ L \}$.
The case (b) does not occur and the stratification $S_2$
of the case (a) consists of only one stratum
$$
 \{ (\aalpha,\bbeta,\ggamma) \,|\, \ggamma \notin \mathbb{Z}_{\leq 0} \}
$$
and the restriction module is isomorphic to $\mathbb{C}$,
because ${\rm ord}_{(-w,w)}(L)=1$ and then $G''=\emptyset$.

Before illustrating steps of the procedure 
$\mbox{\tt Rfr}(\Hg,(\aalpha,\bbeta,\ggamma),G,\ggamma \in \mathbb{Z}_{\leq 0})$,
we show a conclusion that is a list of the restrictions depending on $\ggamma$. 
\begin{enumerate}
\item[(1)]
When $\ggamma \notin \mathbb{Z}_{\leq 0}$,
the restriction module is isomorphic to $\mathbb{C}$.

\item[(2)]
Suppose that $\ggamma = 0$. 
We have $s_0 = 1$ and then $\mathcal{B}_1 = \{1, \p_x\}$ 
(see, e.g., \cite[Alg. 5.2.8]{SST}).
Consider the $\mathbb{C}$-vector space with a basis $\mathcal{B}_1$
that is
$\mathbb{C}^2 = \mathbb{C} \cdot 1 + \mathbb{C} \cdot \p_x$.
Sorting the terms in $L(\aalpha, \bbeta, 0)$ by $<_{(-1,1)}$,
we have
$$
L(\aalpha, \bbeta, 0) = 
x \p_x^2 - x^2 \p_x^2 + (-\aalpha - \bbeta - 1 ) x \p_x - \aalpha \bbeta.
$$
Since the $(-1,1)$-degree of it is $1$, the vector space 
of the denominator of the restriction module is generated by
$(L(\aalpha, \bbeta, 0))|_{x=0}$ 
that is
$$
(L(\aalpha, \bbeta, 0))|_{x=0} = -\aalpha \bbeta .
$$
Thus, the restriction module is 
$\mathbb{C}^2/V$ where 
$
V = \mathbb{C} \cdot \aalpha \bbeta 
$.
Therefore, we have two cases as
  \begin{itemize}
   \item
   When $\aalpha = 0$ or $\bbeta = 0$, $V=\{0\}$ and the restriction module is 
a $2$-dimensinal vector space $\mathbb{C}^2$.
   \item
   When $\aalpha \neq 0$ and $\bbeta \neq 0$, $V=\mathbb{C}$  and the restriction module is a $1$-dimensional vector space $\mathbb{C}$.
   \end{itemize}
\item[(3)]
When $\ggamma \in \mathbb{Z}_{< 0}$,
we can reduce cases of 
$\ggamma \in \mathbb{Z}_{< 0}$ to the case of $\ggamma = 0$
by utilizing left $D$-module isomorphism
$$
D/\Hg(\aalpha, \bbeta, \ggamma) \cong D/\Hg(\overline{\aalpha}, \overline{\bbeta}, \ggamma+1)
$$ 
where 
$\overline{\aalpha}, \overline{\bbeta}$ are 
$\aalpha$ or $\aalpha+1$ and $\bbeta$ or $\bbeta+1$ respectively.
We will prove this fact in Proposition \ref{prop:c-isom-gauss-hg}
\end{enumerate}

We will discuss on  
a stratification of 
$\{(\aalpha,\bbeta,\ggamma)\,|, \ggamma \in \mathbb{Z}_{<0}\}$
after the proposition.

\begin{proposition} \label{prop:c-isom-gauss-hg}
When $\ggamma \in \mathbb{Z}_{< 0}$,
$$D/\Hg(\aalpha, \bbeta, \ggamma) \cong D/\Hg(\overline{\aalpha}, \overline{\bbeta}, 0)$$ 
holds
where 
$\overline{\aalpha}, \overline{\bbeta}$ are 
$\aalpha$ or $\aalpha+1$ and $\bbeta$ or $\bbeta+1$ respectively.
\end{proposition}

\begin{proof}
We abbreviate the Gauss hypergeometric operator as
$L(\ggamma) = L(\aalpha, \bbeta, \ggamma)$
and the left ideal generated by $L$ as 
$\Hg(\ggamma) = \Hg(\aalpha, \bbeta, \ggamma)$. 

The down-step operator $B(\ggamma)$ with respect to $\ggamma$
satisfies
$$\exists P \in D ~{\rm s.t.}~ L(\ggamma-1) B(\ggamma) = P L(\ggamma). $$
The operator
$$
B(\ggamma) = \theta_x +  (\ggamma - 1)
$$
satisfies it.
The up-step operator $H(\ggamma)$ satisfies 
$$L(\ggamma+1) H(\ggamma) = P L(\ggamma).$$
The operator
$$
\Hg(\ggamma) =  (1-x) \p_x + (\ggamma-\aalpha-\bbeta)
$$
satisfies it.

Composing left $D$-module homomorphisms
$$
\varphi : D/\Hg(\ggamma+1) \ni [P] \mapsto [P \cdot H(\ggamma)] \in D/\Hg(\ggamma)
$$
$$
\psi : D/\Hg(\ggamma) \ni [P] \mapsto [P \cdot B(\ggamma+1)] \in D/\Hg(\ggamma+1),
$$
we have
$$
\varphi \circ \psi : D/\Hg(\ggamma) \ni [P] \mapsto [P \cdot B(\ggamma+1) \cdot H(\ggamma)] \in D/\Hg(\ggamma)
$$ 
$$
B(\ggamma+1) \cdot H(\ggamma) \equiv (\aalpha-\ggamma)(\bbeta-\ggamma) ~{\rm mod}~  \Hg(\ggamma)
$$
$$
\varphi \circ \psi = (\aalpha-\ggamma)(\bbeta-\ggamma) {\rm id}.
$$
Reversing the order of the composition,
we have
$$
\psi \circ \varphi : D/\Hg(\ggamma+1) \ni [P] \mapsto [P \cdot H(\ggamma) \cdot B(\ggamma+1)] \in D/\Hg(\ggamma+1)
$$
$$
H(\ggamma) \cdot B(\ggamma + 1) \equiv (\aalpha-\ggamma)(\bbeta-\ggamma) ~{\rm mod}~  \Hg(\ggamma+1)
$$
$$
\psi \circ \varphi  = (\aalpha-\ggamma)(\bbeta-\ggamma)  {\rm id}. 
$$
Hence, when 
$(\aalpha - \ggamma)(\bbeta - \ggamma) \neq 0$,
we have the isomorphism 
$
D/\Hg(\aalpha, \bbeta, \ggamma) \cong D/\Hg(\aalpha, \bbeta, \ggamma+1)
$.

The isomorphism breaks when
$$\aalpha - \ggamma = 0 \text{ or } \bbeta - \ggamma = 0.$$
We derive contiguity relations with respect to $\ggamma$ 
for these cases.

\begin{itemize}
\item[(1)] When $\aalpha - \ggamma = 0$,
put $\aalpha = \ggamma$.
The up-step and down-step operators
for $\Hg(\ggamma, \bbeta, \ggamma)$ with respect to $\ggamma$ 
are
\begin{align*}
&B(\ggamma) = (1-\ggamma)(x(x-1) \p_x +\bbeta x - \ggamma + 1), \\
&H(\ggamma) = (x-1) \p_x + \ggamma, \\
&B(\ggamma+1) \cdot H(\ggamma) \equiv \ggamma^2 (\bbeta - \ggamma) ~{\rm mod}~ \Hg(\ggamma+1), \\
&H(\ggamma) \cdot B(\ggamma+1) \equiv \ggamma^2 (\bbeta - \ggamma) ~{\rm mod}~ \Hg(\ggamma).
\end{align*}
Hence, when $\ggamma^2 (\bbeta - \ggamma) \neq 0$,
$D/\Hg(\ggamma, \bbeta, \ggamma) \cong D/\Hg(\ggamma+1, \bbeta, \ggamma+1)$
holds. 

\item[(1-1)] When $\bbeta - \ggamma = 0$,
put $\bbeta = \ggamma$.
The up-step and down-step operators for  $\Hg(\ggamma, \ggamma, \ggamma)$
with respect to $\ggamma$ are
\begin{align*}
&B(\ggamma) = (1-\ggamma)(x(x-1) \p_x + (2\ggamma-1)x - \ggamma+1), \\
&H(\ggamma) = \p_x, \\
&B(\ggamma+1) \cdot H(\ggamma) \equiv \ggamma^3 ~{\rm mod}~ \Hg(\ggamma+1), \\
&H(\ggamma) \cdot B(\ggamma+1) \equiv \ggamma^3  ~{\rm mod}~ \Hg(\ggamma).
\end{align*}
Hence, when $\ggamma^3 \neq 0$,
$D/\Hg(\ggamma, \ggamma, \ggamma) \cong D/\Hg(\ggamma+1, \ggamma+1, \ggamma+1)$
holds.

\item[(2)] When $\bbeta - \ggamma = 0$,
put $\bbeta = \ggamma$.
The up-step and down-step operators for $\Hg(\aalpha, \ggamma, \ggamma)$
with respect to $\ggamma$ are
\begin{align*}
&B(\ggamma) = (1-\ggamma)(x(x-1) \p_x +\aalpha x - \ggamma + 1), \\
&H(\ggamma) = (x-1) \p_x + \ggamma, \\ 
&B(\ggamma+1) \cdot H(\ggamma) \equiv -\ggamma^2 (\aalpha - \ggamma) ~{\rm mod}~ \Hg(\ggamma+1), \\
&H(\ggamma) \cdot B(\ggamma+1) \equiv -\ggamma^2 (\aalpha - \ggamma) ~{\rm mod}~ \Hg(\ggamma).
\end{align*}
When $\ggamma^2 (\aalpha - \ggamma) \neq 0$,
$D/\Hg(\aalpha, \ggamma, \ggamma) \cong D/\Hg(\aalpha, \ggamma+1, \ggamma+1)$
holds. 

\item[(2-1)] The case $\aalpha - \ggamma = 0$ 
is reduced to the case (1-1).
\end{itemize}
\end{proof}

\begin{example}\rm \label{ex:2F1-Rfr}
Let us illustrate the behavior of $\mbox{\tt Rfr}(\Hg(\beta),p=(a,b,c),G,I)$
where $\beta=(a,b,c)$,
$G=\{ L \}$ and the condition $I$ is $c \in \ZZ_{\leq 1}$.
We retain symbol names of Algorithm \ref{alg:hg-rest}
and of the proof of Proposition \ref{prop:c-isom-gauss-hg}.
The roots of the $b$-function for restriction
are $r_1=0$ and $r_2=1-c$.
Under the assumption $I$, we have $r_2 \geq r_1$.
Then we put
$$ 
 q_1(p)=1-c, q_2(p)=a, q_3(p)=b. 
$$
We can set $\delta=1$ and then $\Lambda=\{ 0 \}$.
The contiguity relation is 
$(B(c), L(c+1), (a-c)(b-c))$.
The $b$-function of contiguity
$(a-c)(b-c)$ can be written in terms of $q$ as
$b(q)=(q_1+q_2-1)(q_1+q_3-1)$.
If $q_2$ and $q_3$ are generic numbers,
there is only one connected component
of (\ref{eq:segments}).
The first coordinate of it is $(-\infty,\infty)$.
Then, the minimum is $0$ which means
$c=1$ ($q_1=0$).
Hence, the restriction module is $\CC$ 
when $c \in \ZZ_{\leq 1}$ and $(a-c)(b-c)\not=0$.
Isomorphisms are given by the contiguity relation.

Let us run $\mbox{\tt Rfr}$ recursively with fewer parameter degrees of freedom.
Let $b_1(q)$ be $q_1+q_3-1$.
We eliminate $q_3$ by $b_1(q)=0$ ($b=c$) and we call
$$
 \mbox{\tt Rfr}(\Hg(q_2,1-q_1,1-q_1),(q_1,q_2),G,
  q_1 \in \ZZ_{\geq 0} \ \mbox{ and }\ q_1+q_3-1=0).
$$
Note that $(q_2,1-q_1,1-q_1)=(a,c,c)$.
Let us execute this procedure.
The roots of $b$-function for restriction is $0$ and $q_1$.
The $\delta$ is $1$ and $\Lambda=\{ 0 \}$.
As we have seen in the proof of Proposition \ref{prop:c-isom-gauss-hg},
the contiguity relation is
$$
\left(
 (x-1)\pd{}+c,
 (1-c)\left( x(x-1)\pd{} + ax - c +1 \right),
 -c^2(a-c)
\right) 
$$
where $-c^2(a-c)=-(1-q_1)^2(q_1+q_2-1)$.
Assume that $q_2=a:={\bar q}_2$ is a generic number.
The connected component of (\ref{eq:segments})
are
$$
(-\infty,1) \times {\bar q}_2, \quad (1,\infty) \times {\bar q}_2.
$$
Then, $Q=\{ 0, 2 \}$.
When $q_1=0$ ($c=1$), the restriction module for $\Hg(a,1,1)$
is a representative of the isomorphic class consisting of $q_1=\{0\}$
and is $\CC$.
When $q_1=2$ ($c=-1$), the restriction module for $\Hg(a,-1,-1)$
is a representative of the isomorphic class consisting of $q_1=\{2,3,\ldots\}$.
Since $s_0=2$ in this case,
the restriction module is
$$
\frac{\CC + \CC \pd{} + \CC \pd{}^2}
 { \langle L|_{x=0}, :\pd{}L: |_{x=0} \rangle }
= \frac{\CC + \CC \pd{} + \CC \pd{}^2}{\CC \pd{}} 
\simeq \CC^2
$$
where $: \ : $ denotes the normally ordered expression
(see, e.g., \cite[p.3]{SST}).
Since the $b$-function for the restriction is 
$-c^2(a-c)=-(1-q_1)^2 (q_1+q_2-1)$,
we need to call recursively $\mbox{\tt Rfr}$
for each factor.
For example, for the factor $1-q_1$,
we call the procedure for
$\Hg(a,0,0)$.
The $b$-function for this contiguity is $a(a+1)$.
Note that the degree of freedom of the parameters decreases when the recursion depth increases.

We believe that these explain how this process works, so we will skip the rest.
\end{example}

\section{Restriction of Appell $F_1$ System to the Origin} \label{sec:f1}

Applying methods discussed in previous sections,
we obtain the following theorem.

\begin{theorem}
The restrictions of the hypergeometric system for the Appell function $F_1$
to $x=y=0$ are as follows.
They are $\CC$-vector spaces.
\begin{itemize}
\item 
When $\ggamma \notin \mathbb{Z}_{\leq 0}$, it is $\mathbb{C}$.

\item 
When $\bbeta = 0$ and $\bbeta' = 0$,
it is
$(\mathbb{C} \cdot 1 + \mathbb{C} \cdot \p_x + \mathbb{C} \cdot \p_y)$.

\item 
When ($\bbeta \neq 0$ or $\bbeta' \neq 0$) and $\aalpha = 0$,
it is
$(\mathbb{C} \cdot 1 + \mathbb{C} \cdot \p_x + \mathbb{C} \cdot \p_y)
/(\mathbb{C} \cdot (-\bbeta' \p_x + \bbeta \p_y)$ .

\item 
When ($\bbeta \neq 0$ or $\bbeta' \neq 0$) and $\aalpha \neq 0$,
it is
$(\mathbb{C} \cdot \p_x + \mathbb{C} \cdot \p_y)/
(\mathbb{C} \cdot (-\bbeta' \p_x + \bbeta \p_y))$.
\end{itemize}
\end{theorem}

Our proof is analogous to the case of the Gauss hypergeometric system.
Several contiguity relations are used.
They are obtained by our implementation
of our algorithms.
Our implementation and details of the proof are published 
in the internet\footnote{ 
\url{https://www.openxm.org/HOME/taka/2025/prog-rest},\\
\url{https://github.com/nobuki-takayama/contiguity2}}.
The proof is omitted here.

\end{document}